\UseRawInputEncoding
\documentclass[preprint,12pt]{elsarticle}
\usepackage{graphicx}
\usepackage{amssymb}
\usepackage{epstopdf}
\usepackage{listings}
\usepackage{float}
\usepackage{lineno}
\usepackage{etex} 
\usepackage{graphicx,color,colordvi}
\usepackage{latexsym,amsbsy,epsfig,fancybox}
\usepackage{amstext}
\usepackage{subfigure}
\usepackage{amsfonts,textgreek}
\usepackage{mathrsfs}
\usepackage{mathtools}
\usepackage{stmaryrd,dsfont}
\usepackage{bbm}
\usepackage{url}
\usepackage{wrapfig,bm}
\usepackage{tikz}
\usepackage{pstricks}
\usepackage{caption}
\usepackage{psfrag}
\usepackage{upgreek}
\usepackage{isomath}
\usepackage{latexsym}
\usepackage{array}
\usepackage{multirow}
\usepackage{booktabs}
\usepackage{verbatim}
\usepackage{color}
\usepackage{colortbl}
\usepackage[colorlinks=true,citecolor=blue]{hyperref}
\usepackage{overpic}
\usepackage{enumerate}
\usepackage{amsmath,amssymb,xspace}
\usetikzlibrary{plotmarks}
\usetikzlibrary{fit,positioning,arrows}
\usetikzlibrary{shapes.geometric, plotmarks}
\usepackage{appendix}

\usepackage[margin=3cm]{geometry}
\usepackage{soul}

\newcommand{\ud}{{\rm d}}
\newcommand{\ui}{{\rm i}}

\tikzstyle{nicebox}=[draw=black!100, fill=white!10, rectangle, inner sep=4pt, inner ysep=16pt]
\tikzstyle{niceboxtitle}=[draw=black!100, fill=white, text=black, rectangle]

\usetikzlibrary{arrows,decorations.markings}
\journal{MSMSE}
\begin{document}
\begin{frontmatter}
\title{Computational Approaches to Model X-ray Photon Correlation Spectroscopy from Molecular Dynamics}
\author[1]{Shaswat Mohanty} 
\author[2]{Christopher B. Cooper}
\author[3]{Hui Wang}
\author[3]{Mengning Liang}
\author[1]{Wei Cai\corref{cor1}}
\ead{caiwei@stanford.edu}
\cortext[cor1]{Corresponding author}
\address[1]{Department of Mechanical Engineering, Stanford University, CA 94305-4040, USA}
\address[2]{Department of Chemical Engineering, Stanford University, CA 94305-4040, USA}
\address[3]{Stanford Linear Accelerator Center, Menlo Park, CA 94025-7015, USA}

\begin{abstract}
X-ray photon correlation spectroscopy (XPCS) allows for the resolution of  dynamic processes within a material across a wide range of length and time scales. 
X-ray speckle visibility spectroscopy (XSVS) is a related method that uses a single diffraction pattern to probe ultrafast dynamics. 
Interpretation of the XPCS and XSVS data in terms of underlying physical processes is necessary to establish the connection between the macroscopic responses and the microstructural dynamics. 
To aid the interpretation of the XPCS and XSVS data, we present a computational framework to model these experiments by computing the X-ray scattering intensity directly from the atomic positions obtained from molecular dynamics (MD) simulations. 
We compare the efficiency and accuracy of two alternative computational methods:
the direct method computing the intensity at each diffraction vector separately, and a method based on fast Fourier transform that computes the intensities at all diffraction vectors at once.
The computed X-ray speckle patterns capture the density fluctuations over a range of length and time scales and are shown to reproduce the known properties and relations of experimental XPCS and XSVS for liquids.
\end{abstract}

\begin{keyword}
X-ray Photon Correlation Spectroscopy \sep X-ray Speckle Visibility Spectroscopy \sep Molecular Dynamics \sep Structure factor \sep Diffusion \sep Optical contrast
\end{keyword}

\end{frontmatter}

\section{Introduction} 
\label{sec:Intro}
X-ray photon correlation spectroscopy (XPCS) is an experimental technique that probes the dynamics of a material at various length and time scales by using the scattering from coherent X-ray sources. It is derived from the Dynamic Light Scattering (DLS) technique which uses laser sources~\cite{berne2000dynamic,chu2008dynamic,zemb2002neutrons}, but the use of coherent X-rays allow us to access much shorter length scales due to the smaller wavelength of X-rays.
The scattering of coherent X-rays from amorphous or disordered material generates seemingly randomly distributed orientations of intensity, called ``speckles", which undergo temporal fluctuations due to atomic motion. The time auto-correlation of the fluctuating speckles (intensity signal) provide information on the characteristic time scales associated with the dynamic processes such as diffusion and stress-relaxation. 
The versatility of the XPCS method has resulted in its use in the study of the dynamics of many material systems, including elastomers~\cite{ehrburger2012xpcs,schlotter2002dynamics,ehrburger2019anisotropic,cooper2020multivalent,nogales2016x}, polymeric solutions~\cite{leheny2012xpcs,guo2012entanglement,papagiannopoulos2005microrheology} and aerogels~\cite{madsen2010beyond,sandy2007contrast,hernandez2014slow}. 
XPCS probes the dynamic processes across a wide range of time scales, ranging from a few seconds to a few $\mu$s~\cite{pecora2008soft,sutton2008review, sinha2014x,leheny2012xpcs,zhang2018dynamics,sandy2018hard,seeck2015x}, with recent experiments showing time resolution at the nanosecond scale~\cite{seaberg2017nanosecond,sun2018pulse}. 
The extension of XPCS in the ultrafast sub $100$-fs regime has also been demonstrated by~\cite{dixon2003speckle,bandyopadhyay2005speckle,decaro2013x} and is referred to as the X-ray speckle visibility spectroscopy (XSVS). The XSVS method enables the measurement of dynamics by tracking the dependence of the optical contrast of the speckle pattern on the exposure time. 
Unlike XPCS, the XSVS method uses the speckle pattern produced by a single (or a double) pulse over a prescribed exposure (or delay) duration.

The coherence length needed for reliable XPCS measurements can be generated by X-ray free-electron lasers (XFELs), available at the Linac coherent light source (LCLS) at the SLAC National Laboratory, USA, and the European X-ray free-electron laser (EuXFEL) in Hamburg, Germany.
XPCS experiments at these facilities can provide critical dynamic information at nm to $\mu$m length scales and ns to ms time scales. 
With the growing availability of XPCS data for complex dynamic systems, there is an urgent need to develop the missing theoretical framework to establish the connection between speckle fluctuations and molecular events.
Establishing such a connection will greatly aid the design of new materials, such as dynamic polymeric networks, aerogels, and recyclable plastics.

The ease of analyzing the results from molecular dynamics (MD) simulations to study the dynamics at different length and time scales has resulted in recent developments in computational XPCS. 
The timescales accessed by XPCS experiments are usually larger than those typically accessed by MD simulations. However, efforts are being made to experimentally probe the dynamics of material systems at the nanosecond time scales that are accessible by MD simulations~\cite{seaberg2017nanosecond,sun2018pulse}. 
Other simulations techniques, such as coarse-grained MD and Monte Carlo, can also capture the material behavior at longer time scales which makes it possible to model the XPCS experiment, albeit at a significant computational cost. 
Nonetheless, to use simulations to link XPCS measurements to molecular scale events would require rapid computation of XPCS signals from a large amount of atomistic configurations.
A direct method exists that computes the scattered intensity from the atomic position in real space obtained from MD simulations. 
Using such a method, the evaluation of the intensity at multiple scattering vectors, $\mathbf{q}$, to study the optical contrast needs to be done individually, so that the computational expense increases linearly with the number of $\mathbf{q}$ points, or pixels on the speckle pattern. 
In the following, we briefly summarize a few of the recent efforts made at developing a computational XPCS and XSVS model to help interpret experimental data. \citet{perakis2017diffusive} describe the XSVS experiment to study diffusion dynamics in the low-to-high density transition in amorphous ice, whereas  \citet{perakis2018coherent} study the slowing dynamics in water due to caging effects at the short length scales by implementing a real-space method to model the XSVS experiment. 
This approach uses the direct method and computes the XPCS signal under the assumption that the Siegert relation holds true. \citet{bikondoa2020x} discuss a Fourier-based numerical framework to simulate the XPCS experiment, with a coarse approximation of the particle density (either 0 or 1) on a grid.

However, the focus of that work is on the modelling of the XPCS instrumentation and does not contain an in-depth discussion of a computational method to obtain the intensity speckles over the entire detector grid.
We are interested in interpreting the XPCS signal in terms of the underlying dynamics across different length scales. Material systems that exhibit both large-scale and small-scale density fluctuations include glass transition phenomena~\cite{berthier2005direct,biroli2004diverging}, hydrophobic phenomena~\cite{vaikuntanathan2014putting,lum1999hydrophobicity}, polymer systems~\cite{zhou2020bridging,patkowski2000long} and dislocation mobility~\cite{groma1997link} to name a few. This calls for an efficient computational XPCS model to interpret the dynamics over multiple length scales, which can be captured by MD simulations involving a large number of atoms.

In this paper, we discuss the comparison between the fast Fourier transform-based (FFT-based) approach and the direct approach to compute the X-ray speckles from atomic configurations to be used in computational XPCS and XSVS models. The FFT-based method is more efficient than the direct approach by simultaneously calculating the scattering intensity over the entire FFT grid. 
By default, our computational XPCS model corresponds to perfectly coherent (in space) and extremely narrow (in time) X-ray pulses, while partial coherency and finite pulse duration can be accounted for through ensemble averaging or time integration.
Through our study, we establish the reliability and efficiency of the FFT-based method in computing the intensity speckles over the entire FFT grid, simultaneously. 
We test the method on liquid Argon (Ar) configurations generated from a molecular dynamics (MD) simulation.
To show that the small size of the simulation box does not produce unwanted artifacts,
we show that the computational XPCS/XSVS data satisfy the known properties/relations that are observed by experiments. Furthermore, we explore the extension of the FFT-based method to estimate the dynamics of water to demonstrate the success of this model in capturing the dynamical information of more complex liquids, previously captured by the direct method.
%

The paper is structured in the following manner. In Section~\ref{sec:Math_form} we describe the mathematical foundation of the scattering theory that is needed for understanding the XPCS and XSVS methods. In Section~\ref{sec:new_method} we discuss the direct computational approach from the atomic positions, along with the algorithm for its implementation. We then present the algorithm for the FFT-based method and discuss its mathematical equivalence to the direct method. 
In Section~\ref{sec:Numerical} we present the computed X-ray speckles from our MD simulations for liquid Ar. We first establish the numerical convergence of the FFT-based method to the direct method, and then and provide numerical evidence for the equivalence of optical contrast in $\bm{q}$-space and in time, and the Siegert relation.  Next, we demonstrate that when considering tracer diffusion of labeled atoms, the time correlation function of speckles indeed decays exponentially with time at a rate proportional to the diffusion constant. Lastly, we show that the optical contrast of the computed X-ray intensities decreases with exposure time, as expected for XSVS experiments. 
As a final benchmark, the FFT-based method is used to compute the XPCS signal of water and successfully reproduces the recent computational results using a different method.
We summarize our discussions in Section~\ref{sec:conclusion}.

\section{Basics of X-ray Scattering and XPCS Theory}  \label{sec:Math_form}
\subsection{Basics of diffraction}
\begin{figure}[H]
    \centering
    {\includegraphics[width=\textwidth]{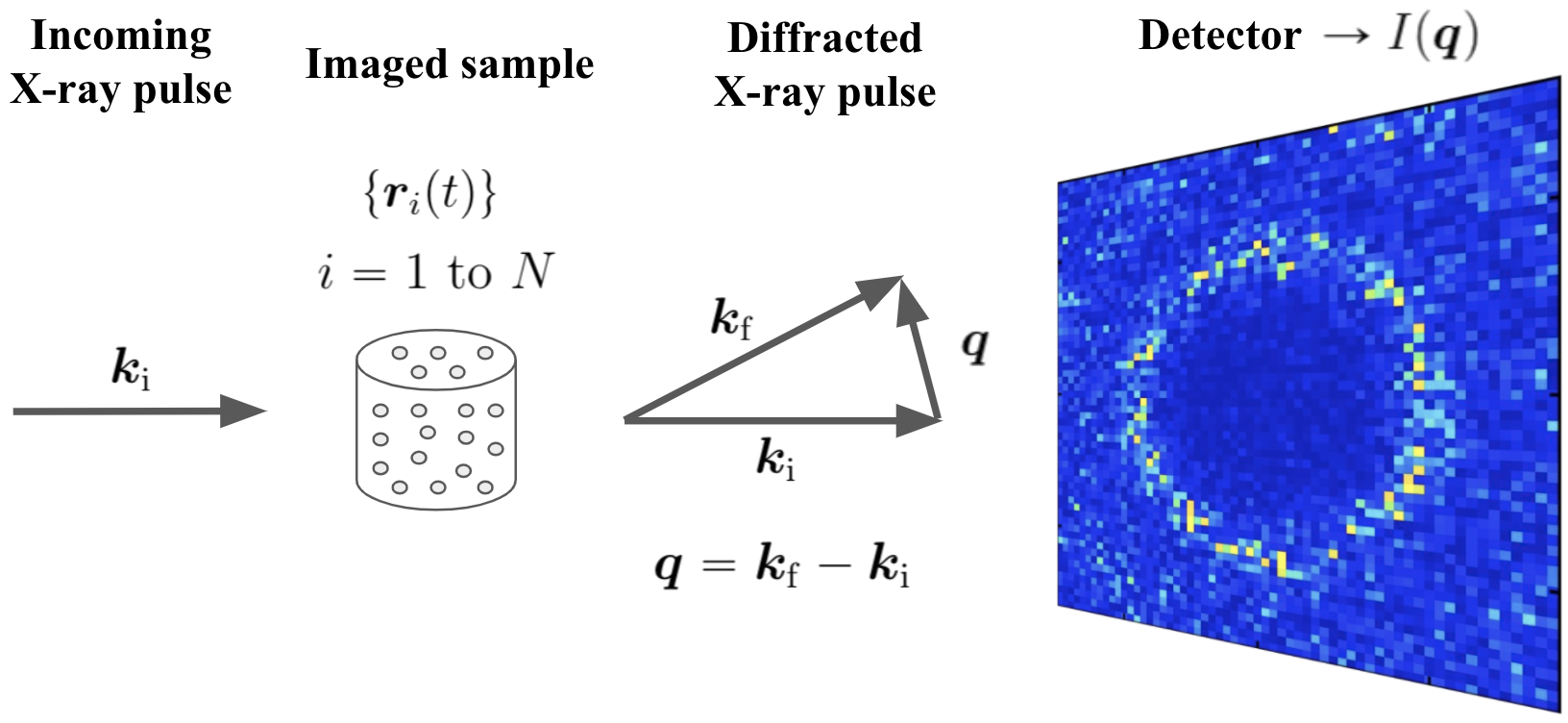}\label{fig:schA}}
    \caption{Schematic of X-ray diffraction in an XPCS experiment. 
    } \label{fig:XPCS_schem}
\end{figure}

In this section, we discuss the key mathematical relations of X-ray scattering theory which provides the foundation for the numerical methods of computational XPCS.
When a coherent X-ray beam with an incoming wavevector, $\bm{k}_{\rm i}$, gets diffracted by the sample with outgoing wavevector $\bm{k}_{\rm f}$, the diffracted wave can be captured on a detector, where the intensity on the detector, $I(\bm{q},t)$, is a function of time $t$ and the change of the wavevector, $\bm{q}=\bm{k}_{\rm f}-\bm{k}_{\rm i}$, as shown in Fig.~\ref{fig:XPCS_schem}. We consider elastic scattering of the incident X-rays from the sample, where the magnitude of the wavevector remain unchanged, and is determined by the X-ray wavelength, $\lambda$,
 \begin{equation}
     \lvert \bm{k}_{\rm i} \rvert = \lvert \bm{k}_{\rm f} \rvert = \frac{2\pi}{\lambda}.
 \end{equation}
The sample consists of a collection of $N$ atoms whose positions as a function of time are denoted by $\{{\bm{r}_{i}(t)}\}$, with $i = 1, 2, \cdots, N$.

In XPCS experiments, the incident X-ray beam consists of very short pulses, each illuminating the sample for a very brief period of time (e.g. $10-120$ fs). 
As an idealization, let us consider the scenario in which the X-ray pulse (at time $t$) is so short that the sample atoms do not move appreciably during this period.
Then the diffracted beam intensity $I(\bm{q},t)$ is determined by the instantaneous atomic positions $\{ \bm{r}_{i}(t) \}$, as given by the following expression \cite{chandler1987introduction,warren1990x},
\begin{equation}
\label{eq:Iq_rij}
    I(\bm{q},t) = \sum_{i=1}^{N}\sum_{j=1}^{N}f_{i}(\bm{q})f_{j}(\bm{q}) \, {\rm e}^{-\ui\bm{q}\cdot[\bm{r}_{i}(t)-\bm{r}_{j}(t)]},
\end{equation}
where $f_{i}(\bm{q})$ is the X-ray atomic form factor of atom $i$. The X-ray atomic form factor can be obtained by the Fourier transform of the electron density field for a given type of atom.  
For many elements, the atomic form factor can be well parameterized by a sum of Gaussians \cite{prince2004international}.

\subsection{Structure factor and scattering intensity}
While Eq.~(\ref{eq:Iq_rij}) already provides the fundamental connection between atomic positions and X-ray diffraction intensities, several additional functions have been introduced in the literature to further clarify the mathematical nature of this connection, which may also lead to more efficient computational methods.
First, let us define the atomic density field, $\rho(\bm{r},t)$, of the sample at time $t$, as a superposition of Dirac delta functions centered at each atom,
\begin{equation}
\label{eq:atom_density}
\rho(\bm{r},t)=\sum_{i}\delta(\bm{r}-\bm{r}_{i}(t))  .  
\end{equation}
In here as well as below, the sum is from $1$ to $N$ unless otherwise specified.
Given the density field $\rho(\bm{r},t)$, we can define its Fourier transform, $p(\bm{q},t)$, as
\begin{align}
    p(\bm{q},t)&= 
    \int \rho(\bm{r},t)\, \, {\rm e}^{-\ui\bm{q\cdot r}}\,\ud^{3}\bm{r} =\sum_{i}\,{\rm e}^{-\ui\bm{q\cdot r}_{i}(t)} .
\end{align}
As we shall see below, $\rho(\bm{r},t)$ and $p(\bm{q},t)$ are only useful in the discussion of X-ray scattering when all atoms in the sample are of the same type.

The inverse Fourier transform of $f_{i}(\bm{q})$ is a measure of the electron density distribution around atom $i$,
\begin{equation}
    \rho^{\rm e}_i(\bm{r}) = \frac{1}{(2\pi)^3}
        \int f_{i}(\bm{q}) \, 
        {\rm e}^{\ui\bm{q}\cdot \bm{r}} \, \ud^3 \bm{q}.
\end{equation}
When the sample contains several types of atoms, it is more useful to consider the electron density with contribution from each atom,
\begin{equation}
\rho^{\rm e}(\bm{r},t)=\sum_{i} 
   \rho^{\rm e}_i(\bm{r}-\bm{r}_{i}(t))  .  
\end{equation}
This can be considered as the field that interacts with the X-ray and generates the diffraction speckles.
Given the electron density field $\rho^{\rm e}(\bm{r},t)$, we can define its Fourier transform, $p^{\rm e}(\bm{q},t)$, as
\begin{equation}
   \label{eq:pe_ri}
    p^{\rm e}(\bm{q},t)=\int \rho^{\rm e}(\bm{r},t)\,
    {\rm e}^{-\ui\bm{q\cdot r}}\,\ud^{3}\bm{r} =\sum_{i}f_{i}(\bm{q})\,{\rm e}^{-\ui\bm{q\cdot r}_{i}(t)} .
\end{equation}
Notice that if all atoms are of the same type, then $p^{\rm e}(\bm{q},t)$ equals to $p(\bm{q},t)$ times the atomic form factor.
Comparing Eqs.~(\ref{eq:Iq_rij}) and (\ref{eq:pe_ri}), we can see that the scattering intensity can also be expressed as,
\begin{equation}
 \label{eq:Iq_pq}
    I(\bm{q},t) = p^{\rm e}(\bm{q},t) \, p^{\rm e *}(\bm{q},t)
     = \left| p^{\rm e}(\bm{q},t) \right|^2,
\end{equation}
where $^{*}$ denotes complex conjugate. An important structural measure of a collection of atoms is the pair distribution function, $g(\bm{r},t)$, defined by,
\begin{equation}
    g(\bm{r},t)=\frac{1}{N\rho_{0}} \sum_{i} \sum_{\substack{j\\ j\neq i}} \delta (\bm{r}-(\bm{r}_{i}(t)-\bm{r}_{j}(t))),
\end{equation}
where $\rho_{0}$ is the bulk atomic density of the sample. 
Although the convention is to exclude the correlation between the atom and itself (i.e. $i=j$) from the definition of $g(\bm{r},t)$, in the following it is more convenient to introduce an alternative definition, $\tilde{g}(\bm{r},t)$, where such a constraint is removed.
\begin{equation}
    \tilde{g}(\bm{r},t)=\frac{1}{N\rho_{0}} \sum_{i} \sum_{j} \delta (\bm{r}-(\bm{r}_{i}(t)-\bm{r}_{j}(t)))
    = g(\bm{r},t)+\frac{1}{\rho_0}\,\delta(\bm{r})
\end{equation}
When all atoms are of the same type, the Fourier transform of $\rho_0\,\tilde{g}(\bm{r},t)$ is the structure factor, $\tilde{S}(\bm{q},t)$,
\begin{equation}
  \label{eq:sq_def_pure}
    \tilde{S}(\bm{q},t)=\frac{1}{N} \sum_{i}\sum_{j} 
    {\rm e}^{-\ui\bm{q}\cdot[\bm{r}_{i}(t)-\bm{r}_{j}(t)]}
     = \frac{1}{N} \, p(\bm{q},t) \, p^*(\bm{q},t).
\end{equation}
When the sample contains atoms of different types, the definition of the structure factor is generalized to the following
\begin{equation}
 \label{eq:Sq_def}
    S(\bm{q},t)=\frac{1}{\sum_{j}f_{j}(\bm{q})^{2}}\sum_{i}\sum_{j} f_{i}(\bm{q})f_{j}(\bm{q})\, {\rm e}^{-\ui\bm{q}\cdot[\bm{r}_{i}(t)-\bm{r}_{j}(t)]}
    = \frac{p^{\rm e}(\bm{q},t) \, p^{\rm e*}(\bm{q},t)} {\sum_{j}f_{j}(\bm{q})^{2}}. 
\end{equation}
We can see that when all atoms are of the same type, Eq.~(\ref{eq:Sq_def}) reduces to Eq.~(\ref{eq:sq_def_pure}).
Comparing Eqs.~(\ref{eq:Iq_rij}) and (\ref{eq:Sq_def}), we can see that the scattering intensity can also be expressed in terms of the structure factor,
\begin{equation}
  \label{eq:Iq_Sq}
    I(\bm{q},t) = S(\bm{q},t) \cdot \sum_{j}f_{j}(\bm{q})^{2}.
\end{equation}
\subsection{Statistical properties of speckles} \label{subsec:stat_prop}
In the XPCS method, the time-varying speckles that are recorded are used to interpret the macroscopic behavior in terms of the dynamics at microscopic length and time scales. Once we obtain the temporal and spatial variation of the XPCS speckle intensity, its time auto-correlation is computed to probe the dynamics of the system at different length scales. The auto-correlation of the speckle at $\bm{q}$ is denoted by $g_{2}(\bm{q}, \tau)$ and it is non-dimensionalized by the square of its mean at $\tau=0$, as given below for infinitesimally short pulses,
\begin{equation}
  \label{eq:g2_def}
    g_{2}(\bm{q},\tau) =\frac{\langle I(\bm{q},t)\, I(\bm{q},t+\tau) \rangle_{t}}{\langle I (\bm{q},t) \rangle_{t}^{2} }
    = \frac{\langle S(\bm{q},t)\, S(\bm{q},t+\tau) \rangle_{t}}{\langle S (\bm{q},t) \rangle_{t}^{2} },
\end{equation}
where $\langle \cdot \rangle_t$ represents the time averaged value of the enclosed entity.
Experimentally, the X-ray pulses always have a finite duration.  Hence the intensity in Eq.~(\ref{eq:g2_def}) should be replaced by the time-averaged X-ray intensity, $I_{\Delta}(\bm{q})$, over the exposure duration of $\Delta_{t}$,
\begin{equation}
  \label{eq:inte_Iq}
    I_{\Delta}(\bm{q}) = \int_{t}^{t+\Delta_{t}} I(\bm{q},t) \,\ud t .
\end{equation}
In this case,
\begin{equation}
\label{eq:g2_def_pulse}
    g_{2}(\bm{q},\tau) =\frac{\langle I_{\Delta}(\bm{q},t)\, I_{\Delta}(\bm{q},t+\tau) \rangle_{t}}{\langle I_{\Delta}(\bm{q},t) \rangle_{t}^{2} } .
\end{equation}
For a X-ray speckle pattern obtained from a single X-ray pulse,  the optical contrast $\beta({q})$ is defined by the variance of the intensity divided by the square of its mean \cite{madsen2010beyond}.
While $\beta(q)$ is given by the scattering intensity distribution $I(\bm{q})$ from an infinitesimally short pulse, we denote the optical contrast from an X-ray pulse of finite duration ${\Delta}_t$ as $\beta_{\Delta}(q)$ where,
\begin{equation}
\label{eq:bq_Iq}
    \beta_{\Delta}(q)=\frac{\langle I_{\Delta}(\bm{q})^{2} \rangle_{q}-\langle I_{\Delta} (\bm{q})\rangle_{q} ^{2}}{\langle I_{\Delta} (\bm{q})\rangle_{q}^{2}}, 
\end{equation}
where $\langle \cdot \rangle_{q}$ represents the average over all detector pixels that satisfy $q-dq/2\le \lvert \bm{q} \rvert < q+dq/2$, for a small $dq$.
For an ergodic and isotropic system, the distribution of pixel intensity at around a particular $\bm{q}$ is the same in $\bm{q}$-space and $t$, so that 
\begin{equation}
\label{eq:bq_bt}
 \langle I_{\Delta}(\bm{q}) \rangle_{q} \approx \langle I_{\Delta}(\bm{q}) \rangle_{t} , \quad \langle I_{\Delta}(\bm{q})^{2} \rangle_{q} \approx \langle I_{\Delta}(\bm{q})^{2} \rangle_{t} .
\end{equation}
This means that we can also define an optical contrast $\beta_0(\bm{q})$ from the time variation of the intensity at a single $\bm{q}$,  
i.e,
\begin{equation}
\label{eq:bq_g2}
    \beta_{0}(\bm{q})=g_{2}(\bm{q},\tau=0)-1. 
\end{equation}
For an ergodic and isotropic system, $\beta_{\Delta}(q)$ and $\beta_0(\bm{q})$ should be equal.

Another useful metric in the study of the dynamic characteristics of the scattering intensity is the intermediate scattering function, $\hat{F}(\bm{q},\tau)$.
When all atoms are of the same type, the intermediate scattering function can be written as,
\begin{equation}
  \tilde{F}(\bm{q},\tau)=\frac{1}{N} \left\langle \sum_{i}\sum_{j} {\rm e}^{-\ui\bm{q}\cdot[\bm{r}_{i}(t)-\bm{r}_{j}(t+\tau)]} \right\rangle_{t} .
\end{equation}
When the sample contains atoms of different types, the definition of the intermediate scattering function is generalized to the following
\begin{equation}
 \label{eq:fq_def_pos}
    {F}(\bm{q},\tau)=\frac{1}{\sum_{j}f_{j}(\bm{q})^{2}} \left\langle
    \sum_{i}\sum_{j} f_{i}(\bm{q})f_{j}(\bm{q})\,  {\rm e}^{-\ui\bm{q}\cdot[\bm{r}_{i}(t)-\bm{r}_{j}(t+\tau)]} \right\rangle_{t} .
\end{equation}
It can be seen that
\begin{align}
\label{eq:Fq_pq}
    {F}(\bm{q},\tau) &=\frac{ \left\langle
      p^{\rm e}(\bm{q},t) \, p^{\rm e*}(\bm{q}, t+\tau)
      \right\rangle_{t} } {\sum_{j}f_{j}(\bm{q})^{2}}.
\end{align}

\noindent We define $\hat{F}(\bm{q},\tau)$ as the normalized intermediate scattering function, i.e.,
\begin{equation}
\label{eq:norm_Fq}
    \hat{F}(\bm{q},\tau) \equiv \frac{{F}(\bm{q},\tau)}{{F}(\bm{q},0)}, 
\end{equation}
so that $\hat{F}(\bm{q},0) = 1$ by definition.
For stationary and ergodic systems, the Siegert relation~\cite{jakeman1974photon,chapman2006femtosecond} connects the auto-correlation function $g_2(\bm{q},t)$, optical contrast, $\beta_0(\bm{q})$, and the normalized intermediate scattering function, $\hat{F}(\bm{q},\tau)$, as follows
\begin{equation}
\label{eq:g2_from_fq}
g_{2}(\bm{q},\tau) =1 + \beta_{0}(\bm{q}) \cdot \lvert \hat{F}(\bm{q},\tau)\rvert ^{2}. 
\end{equation}
The Siegert relation can be derived based on the assumption that the electric field as a function of time of the scattered light is a Gaussian process~\cite{ferreira2020connecting}.
If the atomic trajectories follow a diffusive process (i.e. Brownian motion), then $\hat{F}(\bm{q},\tau)$ reduces to a single exponential,
\begin{equation}
\label{eq:fq_exp}
 \hat{F}(\bm{q},\tau) = \exp[-\Gamma(\bm{q})\, \tau] , 
\end{equation}
where $\Gamma(\bm{q})=D \, q^2$, $q = |\bm{q}|$ and $D$ is the diffusion coefficient.
When the scattering centers are subjected to more complex interactions, the intermediate scattering function can sometimes be modeled by a stretched exponential, 
\begin{equation}
\label{eq:brown}
F(\bm{q},\tau) = \exp[-(\Gamma(\bm{q})\, \tau)^\gamma] ,  
\end{equation}
so that the time correlation of scattered intensity becomes
\begin{equation}
\label{eq:other}
    g_{2}(\bm{q},\tau) =1 + \beta_{0}(\bm{q}) \exp[-2(\Gamma(\bm{q})\, \tau)^{\gamma}]. 
\end{equation}
In Section~\ref{sec:Numerical}, we will test whether the optical contrasts computed from $\bm{q}$-space, Eq.~(\ref{eq:bq_Iq}), and from $t$, Eq.~(\ref{eq:bq_g2}), are equivalent,
as well as the validity of the Siegert relation, Eq.~(\ref{eq:g2_from_fq}), for computational XPCS.
The Siegert relation can be leveraged to express the optical contrast, $\beta_{\Delta}(q)$, in terms of the $\hat{F}(\bm{q},\tau)$ in the case of XSVS experiments~\cite{dixon2003speckle,bandyopadhyay2005speckle,decaro2013x},
\begin{equation}
\label{eq:beta_alt_int}
    \beta_{\Delta}(q) = 2 \beta_{0}(q)\int_{0}^{\Delta_{t}}\left(1-\frac{t}{\Delta_{t}}\right)\lvert \hat{F}(\bm{q},\tau)\rvert ^{2} \frac{\rm dt}{\Delta_{t}}. 
\end{equation}
In Section~\ref{sec:Numerical}, we will test whether the optical contrast computed from Eqs.~(\ref{eq:bq_Iq}), without using the Siegert relation, is equivalent to that computed from Eq.~(\ref{eq:beta_alt_int}) using the Siegert relation.
\section{Computational Approach}
\label{sec:new_method}
From Section~\ref{sec:Math_form}, it can be seen that there are multiple ways to compute the scattered intensity $I(\bm{q},t)$ (e.g. Eqs.(\ref{eq:Iq_rij})),(\ref{eq:Iq_pq}),(\ref{eq:Iq_Sq})) and the auto-correlation function $g_2(\bm{q},t)$ (e.g. Eq.~(\ref{eq:g2_def})).
In the following, we will mostly focus on methods that first compute the scattering intensity $I(\bm{q},t)$ and then compute its auto-correlation function $g_2(\bm{q},t)$, analogous to the workflow of an actual XPCS experiment.
In Section~\ref{subsec:traditional}, we describe the direct method, which computes the structural factor $S(\bm{q},t)$ directly from the atomic positions.
In Section~\ref{subsec:proposed}, we describe a fast Fourier transform (FFT)-based method, which evaluates the electron density field $\rho^{\rm e}(\bm{r},t)$ and its Fourier transform $p^{\rm e}(\bm{q},t)$.
The FFT-based method can also be used to efficiently compute the intermediate scattering function, if needed, using Eq.~(\ref{eq:Fq_pq}).
For simplicity, in this section we shall assume that all atoms are of the same type, so that $f_i(\bm{q}) = f(\bm{q})$ for all $i$.  It is straightforward to generalize the methods described here to samples containing different types of atoms.

\subsection{Direct method} \label{subsec:traditional}
A natural approach to compute the scattering intensity is to first compute the structure factor $\tilde{S}(\bm{q},t)$, defined in Eq.~(\ref{eq:sq_def_pure}), and then multiply the atomic form factors, as in Eq.~(\ref{eq:Iq_Sq}).
However, a brute-force implementation of Eq.~(\ref{eq:sq_def_pure}) from atomic positions is usually infeasible because it requires a summation over all atomic pairs, with a computational cost of $\mathcal{O}(N^{2})$.
Furthermore, MD simulations usually assume periodic boundary conditions; accounting for all the periodic images for each atom leads to an infinite sum that is usually intractable unless truncated.

To compute the intensity at any arbitrary wavevector, $\bm{q}$, it is customary to introduce a cut-off radius, $r_{\rm cut}$, so that only atomic pairs $(i, j)$ with $\lvert \bm{r}_{i}-\bm{r}_{j} \rvert < r_{\rm cut}$ are explicitly included in the sum, while contributions from atomic pairs of greater distances are accounted for with an approximation.
Such an approach is indeed often used in the calculation of the pair distribution function $\tilde{g}(\mathbf{r})$, whose Fourier transform gives $S(\mathbf{q})$.
This approach usually gives results with acceptable accuracy for the angularly averaged structural factor, $S(q)$, which varies smoothly over the magnitude $q$ of the $\mathbf{q}$-vector.
However, we have found that introducing a (sharp) real-space cut-off radius, $r_{\rm cut}$, produces unacceptably large (aliasing) error in $S(\mathbf{q})$, which varies strongly from one pixel to the next in $\mathbf{q}$-space.
Alternatively, the structure factor can be obtained from the $p(\mathbf{q},t)$ and its complex conjugate, where $p(\mathbf{q},t)$ can be calculated from a single sum over all atomic positions,
\begin{equation}
\label{eq:sq_cut_simple}
    \tilde{S}(\bm{q},t)=\frac{1}{N}\underbrace{\sum_{i}{\rm e}^{-\ui\bm{q}\cdot\bm{r}_{i}(t)}}_{p(\bm{q},t)}\underbrace{\sum_{j}{\rm e}^{\ui\bm{q}\cdot\bm{r}_{j}(t)}}_{p^*(\bm{q},t)} .
\end{equation}
The periodic boundary condition of the simulation cell is naturally taken care of as long as each atom is summed only once, and if $\mathbf{q}$ lies on the reciprocal lattice of the Bravais lattice defined by the three repeat vectors of the simulation box.
In the simplest case of a cubic simulation box of length $L$, the requirement for $\mathbf{q}$ is that each of its component must be an integer multiple of $2\pi/L$.
In this case, it does not matter which one of the periodic images of any atom is used in the summation for computing $p(\mathbf{q},t)$.
A consequence of the periodic boundary condition is that the scattered intensity is strictly zero if $\bm{q}$ does not lie exactly on the reciprocal lattice.
The direct method can be summarized by the following algorithm.

\vspace{0.1in}

\noindent {\it Algorithm 1} 
\begin{itemize}
\item Construct a grid of the discretized wavevectors, $\bm{q}$, which are the reciprocal lattice of the Bravais lattice defined by the repeat vectors of simulation box. (For a cubic simulation box with length $L$, $\mathbf{q}$ is a simple cubic lattice whose components are integer multiples of ${2\pi}/{L}$.) 
\item Compute the $p(\bm{q})$ from the atomic positions, by summing ${\rm e}^{-i\bm{q}\cdot\bm{r}_{i}}$ over each atom $i$.
\item Compute the structural factor at the chosen $\bm{q}$ using Eq.~(\ref{eq:sq_cut_simple}), and compute the X-ray intensity using Eq.~(\ref{eq:Iq_Sq}).
\end{itemize}
An advantage of {\it Algorithm 1} is that the computations of scattered intensity for different $\bm{q}$ vectors are independent of each other, and hence can be performed in parallel.
However, this is also related to a disadvantage of this method --- the total computational cost increases linearly with the number of $\bm{q}$ vectors of interest.
Because each pixel in the X-ray detector covers a small but finite angular range of the diffracted X-ray, an integral over the $\bm{q}$ space may be needed even to simulate the scattered X-ray signal detected by a single pixel.
Performing such an integral numerically requires evaluation of $I(\bm{q},t)$ at multiple $\bm{q}$ vectors and can significantly increase the computational cost.
\subsection{FFT-based method} \label{subsec:proposed}

For the FFT-based approach, we compute a modified atomic density ${\rho}^{\eta}(\bm{r},t)$, whose Fourier transform can be used to obtain the Fourier transform of the electron density field,  ${p}^{\rm e}(\bm{q},t)$, so that we can  obtain the scattering intensity $I(\bm{q},t)$ from Eq.~(\ref{eq:Iq_pq}). Theoretically, the atomic density is given by Eq.~(\ref{eq:atom_density}), but resolving the these $\delta$-function peaks on a numerical grid is not feasible. 
Therefore, we consider a ``smeared-out'' atomic density field, as a superposition of set of a 3-dimensional Gaussian distributions centered at each atomic position. 
\begin{equation}
\label{eq:modified_r}
  {\rho}^\eta(\bm{r},t)=\sum_{i}\left(\frac{1}{\sqrt{2\pi}\eta}\right)^{3} {\rm e}^{-\frac{|\bm{r}-\bm{r}_{i}(t)|^{2}}{2\,\eta^{2}}}, 
\end{equation}
where $\eta$ is the standard deviation of the Gaussian distribution; it is a numerical parameter that should be chosen appropriately (see below).
The resulting density field can be represented on a sufficiently fine grid.
For simplicity, we shall consider a cubic simulation cell (subjected to periodic boundary condition) and representing the density field on an $N_{\rm grid}\times N_{\rm grid}\times N_{\rm grid}$ grid.
Computing the atomic density field can be expensive if a Gaussian function centered on each atom needs to be evaluated on all $N_{\rm grid}^{3}$ grid points. 
However, because a Gaussian distribution decays rapidly with distance, we can limit the evaluations of the Gaussian distribution in a cubic region containing $k\times k\times k$ grid points (accounting for periodic boundary conditions) without introducing much numerical error. To this end we choose $k$ as the nearest integer greater than $10\eta$ to successfully capture the features of the distribution to up to $\pm5\eta$.
This limits the computational cost of constructing the density field ${\rho}^\eta(\bm{r},t)$ to $\mathcal{O}(k^3 N_{\rm grid})$.
The imposition of the truncation of the cubic grid is similar to the Particle Mesh Ewald (PME) method~\cite{darden1999new} in computing Coulomb interactions in atomistic simulations using periodic supercells. In the PME method, the quick convergence of the potential and charge density fields in the real and Fourier space result in a negligible loss of accuracy with the introduction of the truncation.

The Fourier transform of ${\rho}^\eta(\bm{r},t)$ is related to the Fourier transform of the atomic density, $p(\bm{q},t)$,  as follows,
\begin{equation}
\label{eq:p_rq_app}
    {p}^\eta(\bm{q},t) = \int {\rho}^{\eta}(\bm{r},t)\,{\rm e}^{-\ui\bm{q\cdot r}}\,\ud^{3}\bm{r} = f^\eta(\bm{q}) \, p(\bm{q},t), 
\end{equation}
where $f^\eta(\bm{q}) = {\rm e}^{-\frac{|\bm{q}|^{2}\eta^{2}}{2}}$ is the Fourier transform of the Gaussian distribution with standard deviation $\eta$.
Therefore, we can obtain the  Fourier transform of the electron density field, ${p}^{\rm e}(\bm{q},t)$, from ${p}^\eta(\bm{q},t)$,
\begin{equation}
\label{eq:pq_rq_app_alt}
    {p}^{\rm e}(\bm{q},t) = f(\bm{q})\,{p}(\bm{q},t)= \frac{f(\bm{q})}{f^\eta(\bm{q})}\,
    {p}^{\eta}(\bm{q},t).
\end{equation}
From the ${p}^{\rm e}(\bm{q},t)$ we can compute the $F(\bm{q},\tau)$ by,
\begin{align}
      \label{eq:Fq_pq_alt}
    {F}(\bm{q},\tau) &=\frac{ \left\langle
      p^{\rm e}(\bm{q},t) \, p^{\rm e*}(\bm{q}, t+\tau)
      \right\rangle_{t} } {\sum_{j}f_{j}(\bm{q})^{2}} = \frac{ 1 } {N\, f^\eta(\bm{q})^2 }\left\langle
        \,{p}^{\eta}(\bm{q},t) \, \,{p}^{\eta *}(\bm{q},t+\tau)\right\rangle_{t}.
\end{align}
The structural factor can now be computed by,
\begin{equation}
    \tilde{S}(\bm{q},t) =\frac{1}{N} \, \,{p}(\bm{q},t) \, {p}^{*}(\bm{q},t)
        =\frac{ 1 } {N\, f^\eta(\bm{q})^2 }
        \,{p}^{\eta}(\bm{q},t) \, \,{p}^{\eta *}(\bm{q},t).
\end{equation}
The scattered X-ray intensity can be computed using,
\begin{equation}
\label{eq:Iq_pq_new}
    I(\bm{q},t) = N \, f(\bm{q})^2 \, \tilde{S}(\bm{q},t)
               = \frac{f(\bm{q})^2} { f^\eta(\bm{q})^2 } 
                \,{p}^{\eta}(\bm{q},t) \, \,{p}^{\eta *}(\bm{q},t). 
\end{equation}
The Fourier transform from ${\rho}^{\eta}(\bm{r},t)$ to ${p}^\eta(\bm{q},t)$ is computed using fast Fourier transform (FFT).
Due to the discretization error of the FFT grid, $f^\eta(\bm{q})$ is not exactly the same as ${\rm e}^{-\frac{|\bm{q}|^{2}\eta^{2}}{2}}$.
We have found that a much higher accuracy is reached if we compute $f^\eta(\bm{q})$ from FFT of a 3d Gaussian distribution centered at the origin with standard deviation $\eta$ (truncated to the nearest $k\times k\times k$ grid points).
Therefore, the fast Fourier transform of the density field ${p}^\eta(\bm{q},t)$ provides an alternative method to compute the intensity of the diffracted X-ray speckles simultaneously on all grid points in the $\bm{q}$-space.
The proposed Fourier-based method can be summarized by the following algorithm, which computes $\tilde{S}(\bm{q})$ and $I(\bm{q})$ given the atomic positions $\{\bm{r}_i\}$ 
and over a 3D grid in the $\bm{q}$-space.
For simplicity, we shall assume that the simulation cell has a cubic shape with length $L$ (subjected to periodic boundary conditions) and is discretized into an $N_{\rm grid}\times N_{\rm grid}\times N_{\rm grid}$ grid.
The user also needs to specify the width of the Gaussian distribution, $\eta$, and the number $k$ of the grid cut-off size around each atom.
\vspace{0.1in}

\noindent {\it Algorithm 2} 
\begin{itemize}
\item Initialize the $\rho^\eta(\bm{r})$ field with zero values on an $N_{\rm grid}\times N_{\rm grid}\times N_{\rm grid}$ grid. 
\item For each atom $i$, compute the Gaussian distribution $\left(\frac{1}{\sqrt{2\pi}\eta}\right)^{3} {\exp}\left({-\frac{|\bm{r}-\bm{r}_{i}|^{2}}{2\,\eta^{2}}}\right)$ on the nearest $k\times k\times k$ grid points (accounting for periodic boundary conditions), and add the values to
the $\rho^\eta(\bm{r})$ field.
\item Compute $p^{\eta}(\bm{q})$ from $\rho^\eta(\bm{r})$ by fast Fourier transform (FFT).  The resulting field is represented on a regular grid in the $\bm{q}$-space. 
\item Compute $f^\eta(\bm{q})$ over the $\bm{q}$-space grid from the FFT of a 3D Gaussian distribution with standard deviation $\eta$, centered at the origin (with zero values outside the nearest $k\times k\times k$ grid points).
\item Compute the scattering intensity over the $\bm{q}$-space grid using Eq.~(\ref{eq:Iq_pq_new}).
\end{itemize}

The entire speckle pattern is generated (so that the optical contrast can be immediately obtained) after executing Algorithm 2, as opposed to at one $\bm{q}$-vector only after executing Algorithm 1. In addition to the $I(\bm{q},t)$, Algorithm 2 can be easily adapted to compute the $F(\bm{q},t)$ from Eq.~(\ref{eq:Fq_pq_alt}), using the ${p}^{\rm e}(\bm{q},t)$ obtained from Eq.~(\ref{eq:pq_rq_app_alt}).

\section{Numerical Results} 
\label{sec:Numerical}

\subsection{Atomistic model}
\begin{figure}[H]
    \centering
    \subfigure[]{\includegraphics[width=0.48\textwidth]{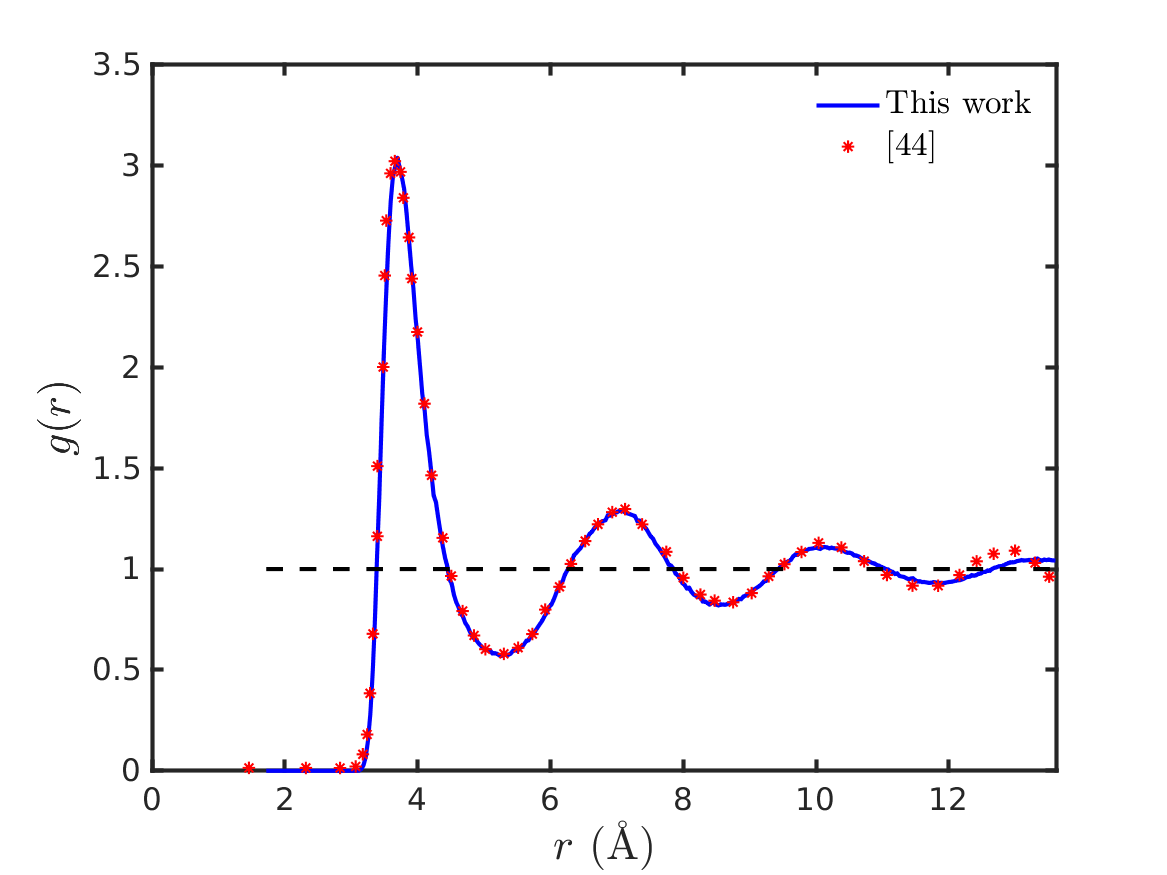}\label{fig:avggrA}} 
    \subfigure[]{\includegraphics[width=0.48\textwidth]{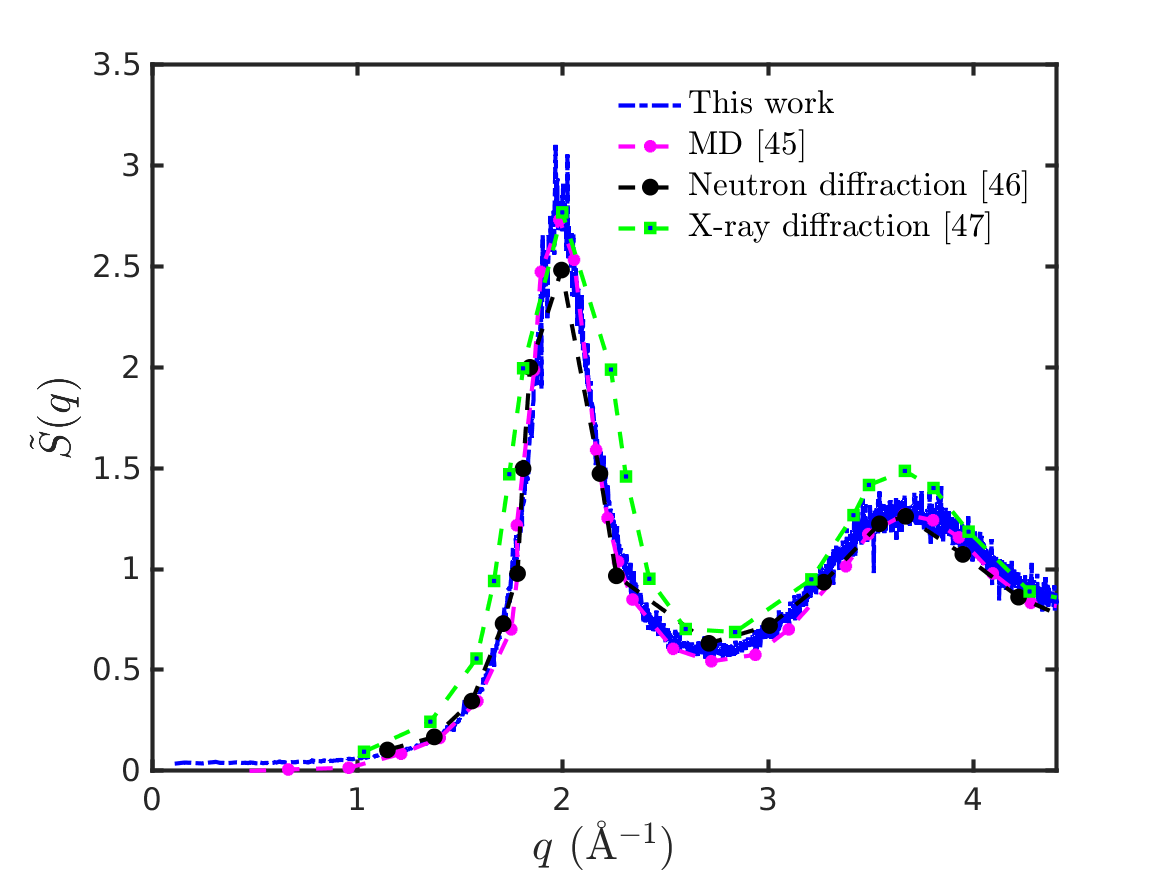}\label{fig:avggrB}} 
    \caption{ The (a) pair correlation function $g(r)$ and (b) structural factor $s(q)$ time-averaged over $2.156$ ns validated against the benchmark results~\cite{rahman1964correlations} from Molecular Dynamics \cite{verlet1968computer}, Neutron diffraction \cite{henshaw1957atomic} and X-ray diffraction \cite{eisenstein1942diffraction}. \label{fig:validation}}
\end{figure}

We choose the simple liquid of Argon (Ar) as a test bed to benchmark the two methods because of readily available previous results from both experiments and computations.
Our MD simulations use the Lennard-Jones potential, which has been shown to provide an accurate description for the interatomic interactions in liquid Ar,
\begin{align}
    V_{\rm LJ}(r) &=4 \epsilon\left[ \left(\frac{\sigma}{r}\right)^{12}-\left( \frac{\sigma}{r}\right)^{6}\right],
\end{align}
where $\epsilon=16.5402 \times 10^{-22}$ J/atom is the depth of the potential well, and $\sigma=3.405$~\AA\, is the characteristic atomic size~\cite{lv2011molecular}.
The simulations are performed at density, $\bar{\rho}=1680.2$ kg/m$^{3}$, and temperature, $T=85.1$ K~\cite{rahman1964correlations}.
The MD simulation box is a cube with length $L = 59.19\,$\AA $\;$ and contains $4000$ atoms. The discretization of the box in real space is given by $\delta=L/N_{\rm grid}$.
The initial atomic positions are randomly positioned, followed by a relaxation to a local energy minimum.
An equilibration MD simulation is then performed using the NPT ensemble for $1.028$ ns. 
The equilibrated configuration is then used as the initial condition for an MD simulation using the NVT ensemble for $2.156$~ns (with a time step of $\Delta t = 10.78$ fs).
Snapshots of atomic positions are recorded every $43.12$ ps.
To validate the computation by the FFT based approach, we compute the $S(q)$ and $g(r)$ from the MD simulation generated trajectory. Fig.~\ref{fig:avggrA} show the time averaged and orientation averaged pair distribution function, $g(r)$, from our MD simulation, which is in good agreement with the literature values~\cite{rahman1964correlations}. These time averages are performed over the last $40$ saved configurations from the MD simulation.
The corresponding time averaged and orientation averaged structural factor, $\tilde{S}(q)$, as shown in Fig.~\ref{fig:avggrB}, also agrees with previous reports from MD simulations~\cite{verlet1968computer}, neutron diffraction~\cite{henshaw1957atomic} and X-ray diffraction~\cite{eisenstein1942diffraction}. 
\subsection{Convergence of FFT-based method to direct method}
While the direct method (Algorithm 1) is essentially exact in the evaluation of scattering intensity at a given $\bm{q}$-vector (except for round-off errors), the FFT-based method (Algorithm 2) introduces several numerical approximations, as described by parameters $\eta$, $N_{\rm grid}$ and $k$.
Here we demonstrate that the errors from these approximations can be made vanishingly small.
\begin{figure}[H]
    \centering
        \subfigure[]{\includegraphics[width=0.48\textwidth]{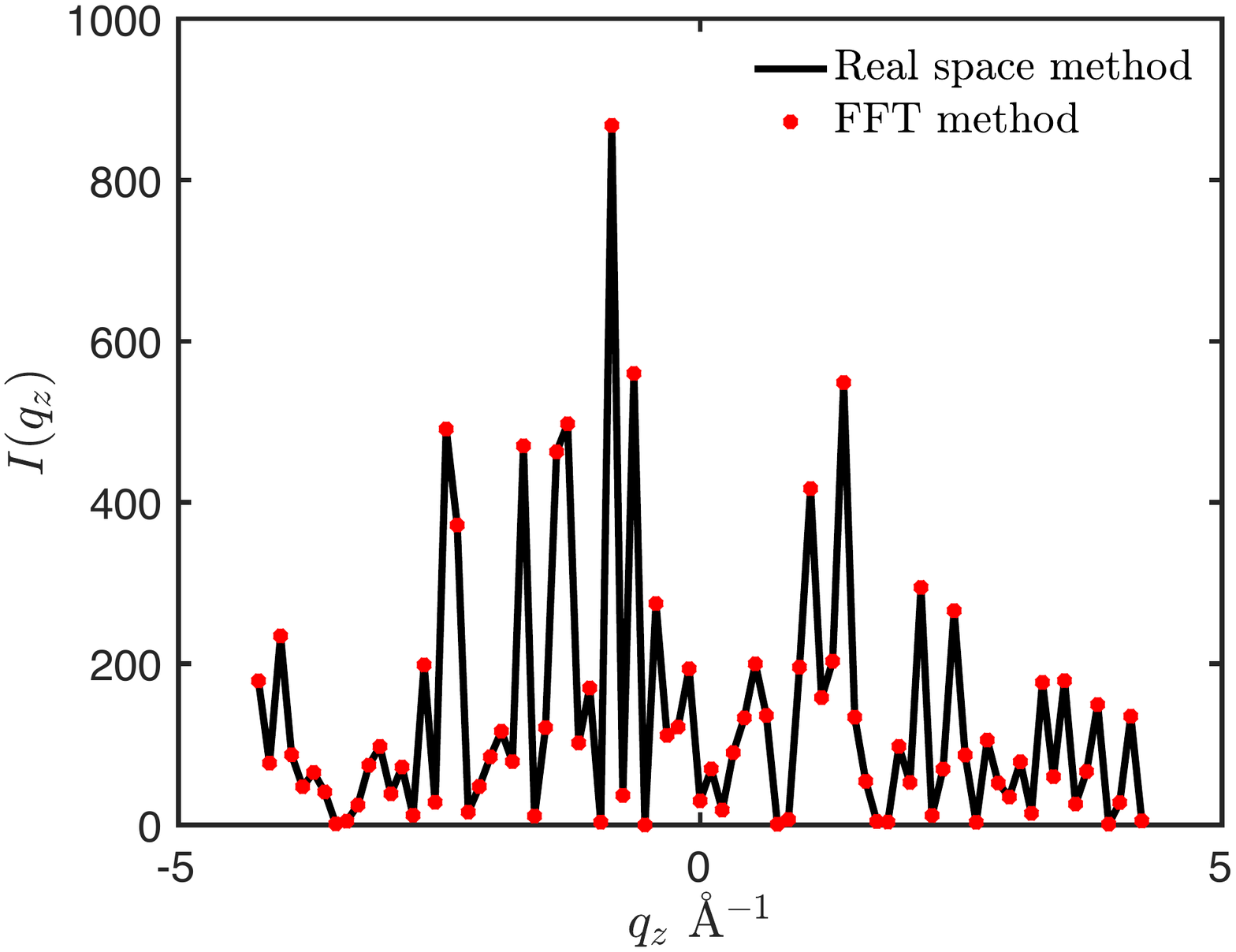}
    }\label{fig:valid_pos}
     \subfigure[]{\includegraphics[width=0.48\textwidth]{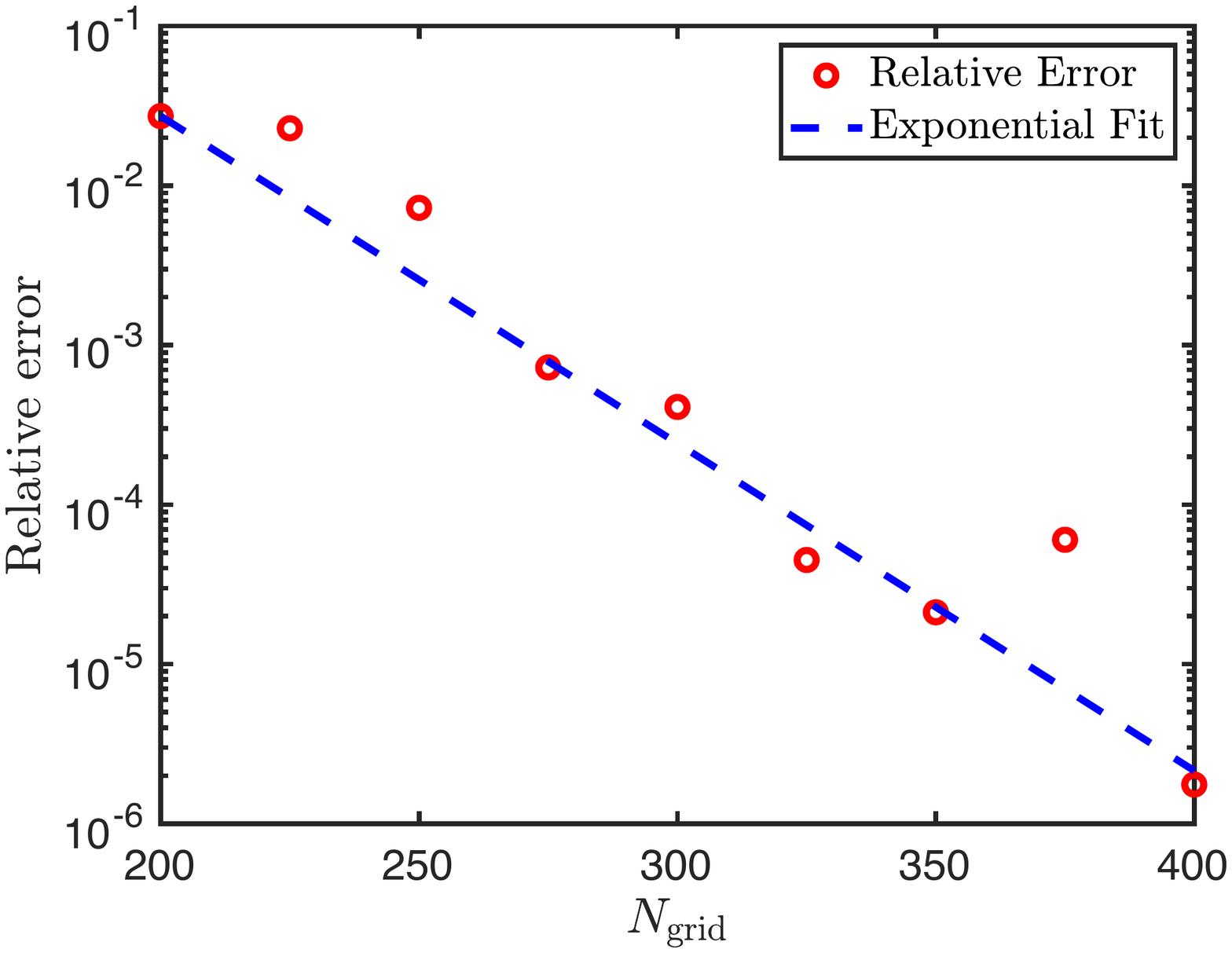}
    \label{fig:conv_pos}}
    \caption{(a) The comparison between the intensity obtained from the real space (direct) and Fourier space method at $\delta = \eta = 0.1485$ \AA, as a function of $q_z$ for $q_x = 0$ \AA$^{-1}$ and $q_y = 1.481$ \AA$^{-1}$. 
    (b) The relative error between the real-space (direct) and FFT-based method as a function of the FFT grid size $N_{\rm grid}$. \label{fig:real_fourier_covnergence}}
\end{figure}

Fig.~\ref{fig:real_fourier_covnergence} show the X-ray intensity (in arbitrary units) computed from the direct method (solid line) and the FFT-based method (dots) as a function of $q_z$ for $q_x = 0$ \AA$^{-1}$ and $q_y = 1.481$ \AA$^{-1}$.
For this numerical example to test the agreement between the two methods, we use $N_{\rm grid} = 400$, $k=11$ and $\eta = 0.1485$ \AA.
Fig.~\ref{fig:real_fourier_covnergence}(a) shows that the results from both methods agree well with each other.
To quantify the error of the FFT-based approach, Fig.~\ref{fig:real_fourier_covnergence}(b) plots the relative difference between the two methods at $q_x = 0$ \AA$^{-1}$, $q_y = 1.481$ \AA$^{-1}$, $q_z = -1.164$ \AA$^{-1}$, as a function of FFT grid size $N_{\rm grid}$.
The error is seen to decay exponentially fast with increasing grid size $N_{\rm grid}$.
We note that in order to achieve this level of convergence, it is important to compute $f^\eta(\bm{q})$ from FFT instead of the analytic expression, as discussed in Section~\ref{subsec:proposed}.

\begin{figure}[H]
    \centering
        \subfigure[]{\includegraphics[width=0.48\textwidth]{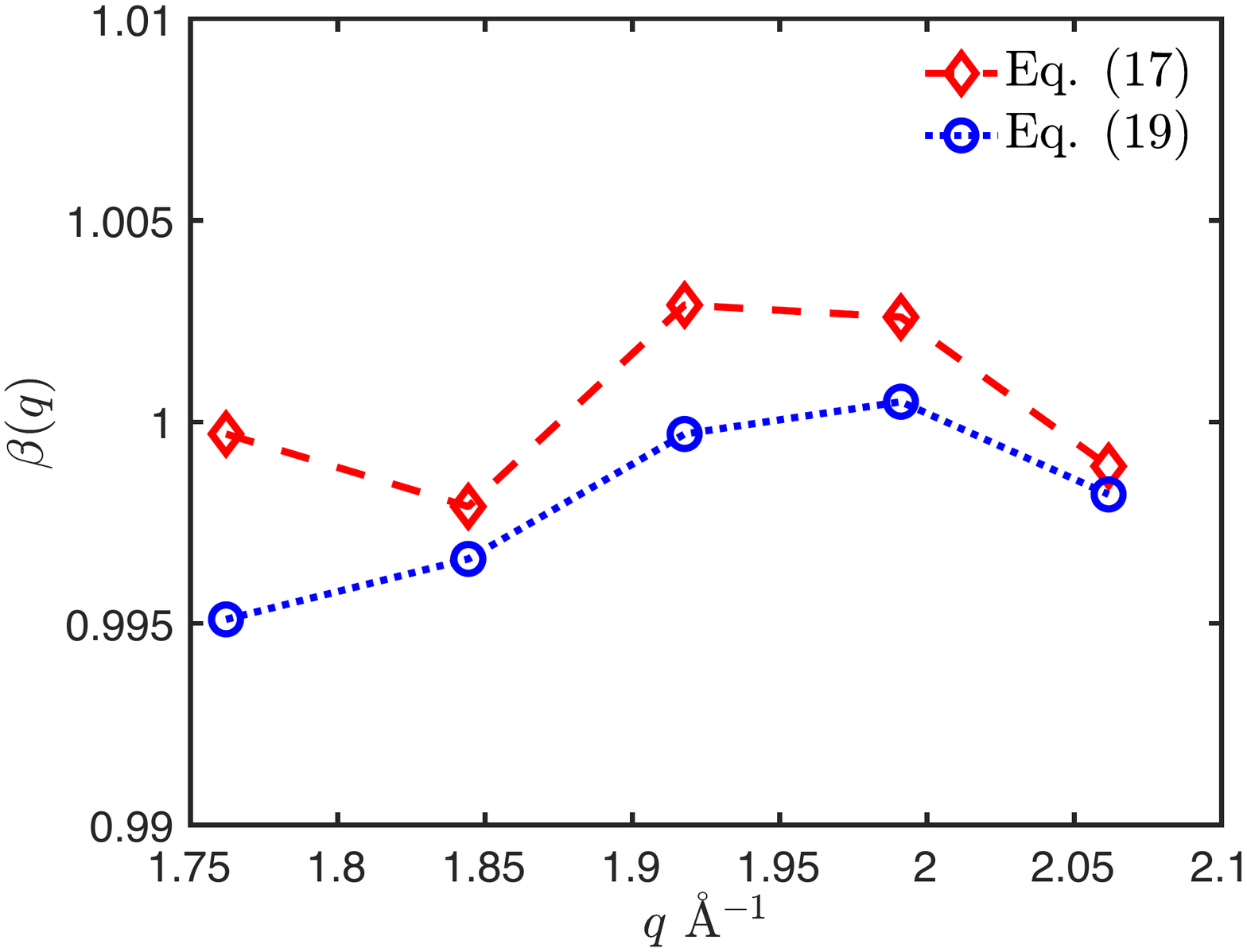}\label{fig:beta_comp}}
        \subfigure[]{\includegraphics[width=0.48\textwidth]{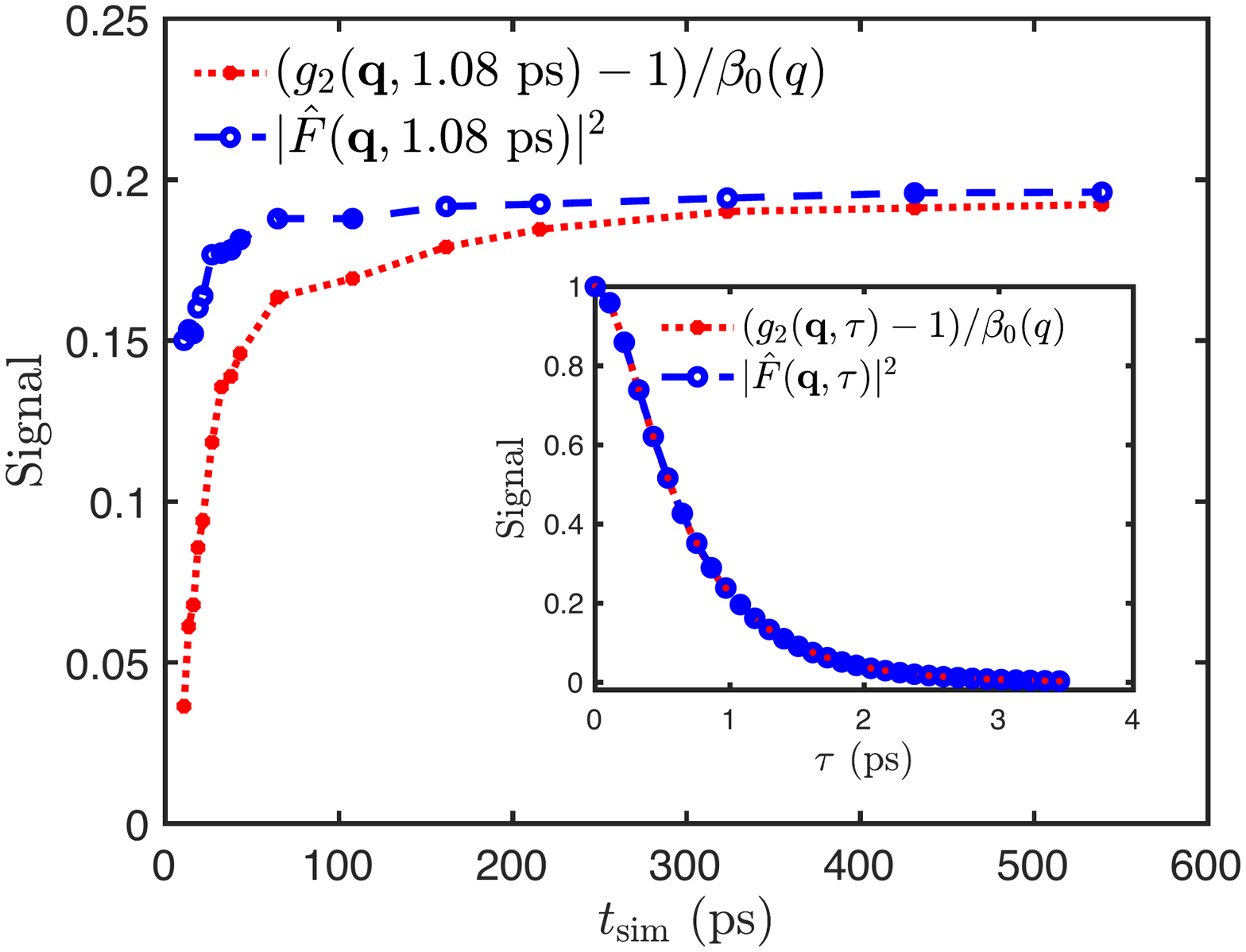}\label{fig:g2_fq}}
    \caption{(a) Comparison between the optical contrast computed from Eqs.~(\ref{eq:bq_Iq}) and (\ref{eq:bq_g2}). (b) Convergence of $(g_{2}(\bm{q},\tau)-1)/\beta_{0}(q)$ and $\lvert \hat{F}(\bm{q}) \rvert ^{2}$ at $\tau = 1.08$ ps with increasing simulation time ($t_{\rm sim}$).
    The inset shows between $(g_{2}(\bm{q},\tau)-1)/\beta_{0}(q)$ and $\lvert \hat{F}(\bm{q},\tau) \rvert ^{2}$ as a function of $\tau$ for $t_{\rm sim} = 539$ ps. All tests are carried out for $\lvert \bm{q} \rvert = 1.844 \pm 0.029$ \AA$^{-1}$.
    }
    \label{fig:siegert}
\end{figure}
To demonstrate the capability of the numerical model to capture the known relations of the XPCS and XSVS experiments, we first test the equivalence of optical contrast defined in Eq.~(\ref{eq:bq_Iq}) and Eq.~(\ref{eq:bq_g2}), and then the Siegert relation, Eq.~(\ref{eq:g2_from_fq}). 
To test this equivalence, we first compute the $\beta(q)$ over the $\bm{q}$-space using Eq.~(\ref{eq:bq_Iq}) for each one of the $5000$ MD snapshots and then averaging the results. 
Second, we compute $\beta(q)$ over time using Eq.~(\ref{eq:bq_g2}) at a given $\bm{q}$ over all time, and then averaging the results over all $\bm{q}$ points 
that fall within the $q -dq/2\le \lvert \bm{q} \rvert < q+dq/2$ ring. 
Fig.~\ref{fig:beta_comp} shows that the results from these two approaches agree well with each other, across different values of $q$, thereby confirming the equivalence between Eqs.~(\ref{eq:bq_Iq}) and (\ref{eq:bq_g2}).
We note that the optical contrast here is very close to 1, which is expected since we assume fully coherent X-rays. This result is a further indication of the successful computation by our numerical model.
To test the Siegert relation, we compare $(g_{2}(\bm{q},\tau)-1 )/ \beta_0(\bm{q})$ from the computed intensity (Eq.~(\ref{eq:g2_def_pulse})), and $\lvert \hat{F}(\bm{q},\tau) \rvert ^{2}$ directly from the atomic positions (Eq.~(\ref{eq:fq_def_pos})). %
Fig.~\ref{fig:g2_fq} shows that the Siegert relation, and hence the stationarity and ergodicity assumptions, are well satisfied for the MD simulation of liquid Ar at $T = 100$~K.

\subsection{Speckle pattern of X-ray diffraction} \label{subsec:prelim_valid}
\begin{figure}[H]
    \centering
    \subfigure[]{\includegraphics[width=0.48\textwidth]{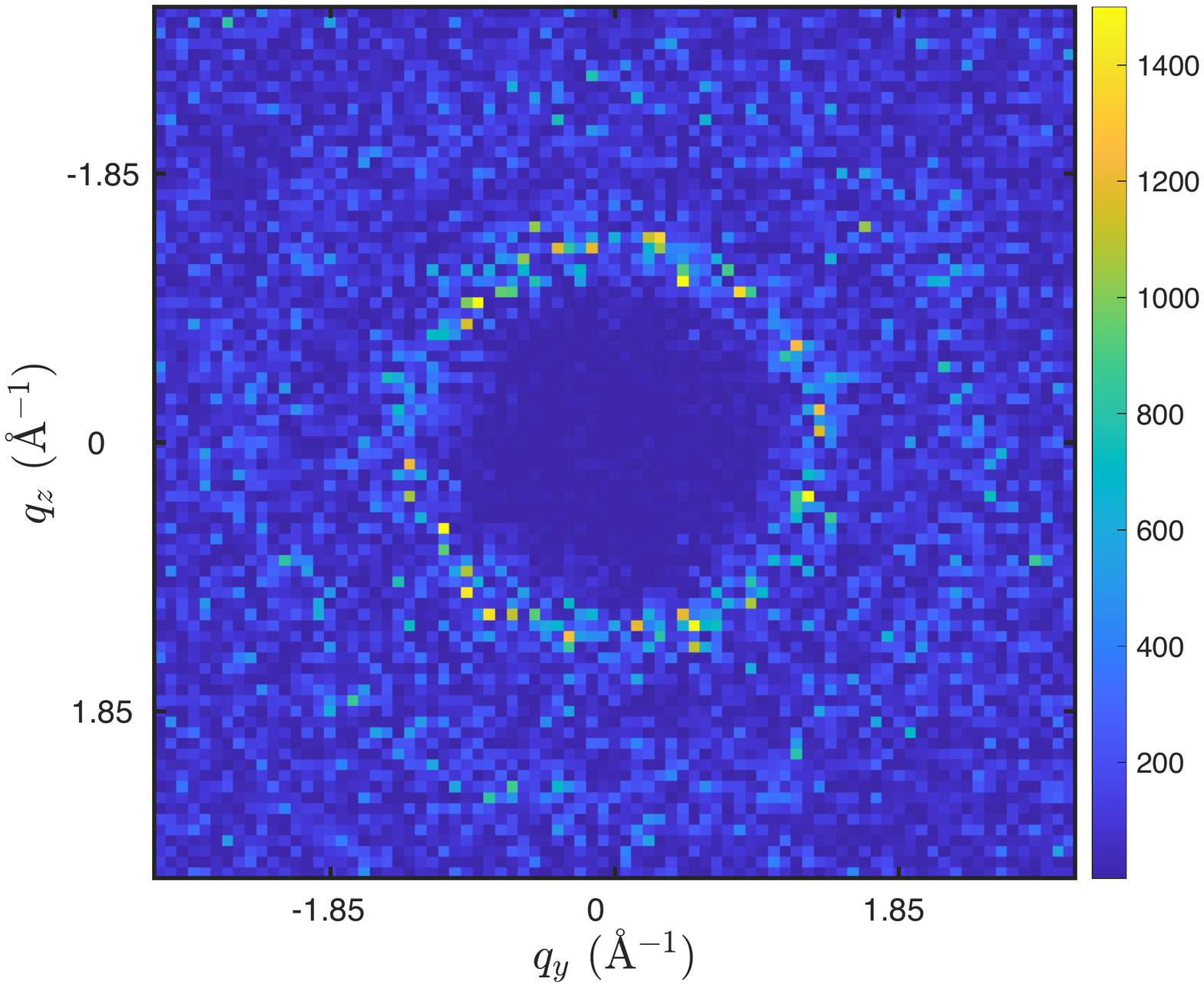}\label{fig:sqA}}
    \subfigure[]{\includegraphics[width=0.48\textwidth]{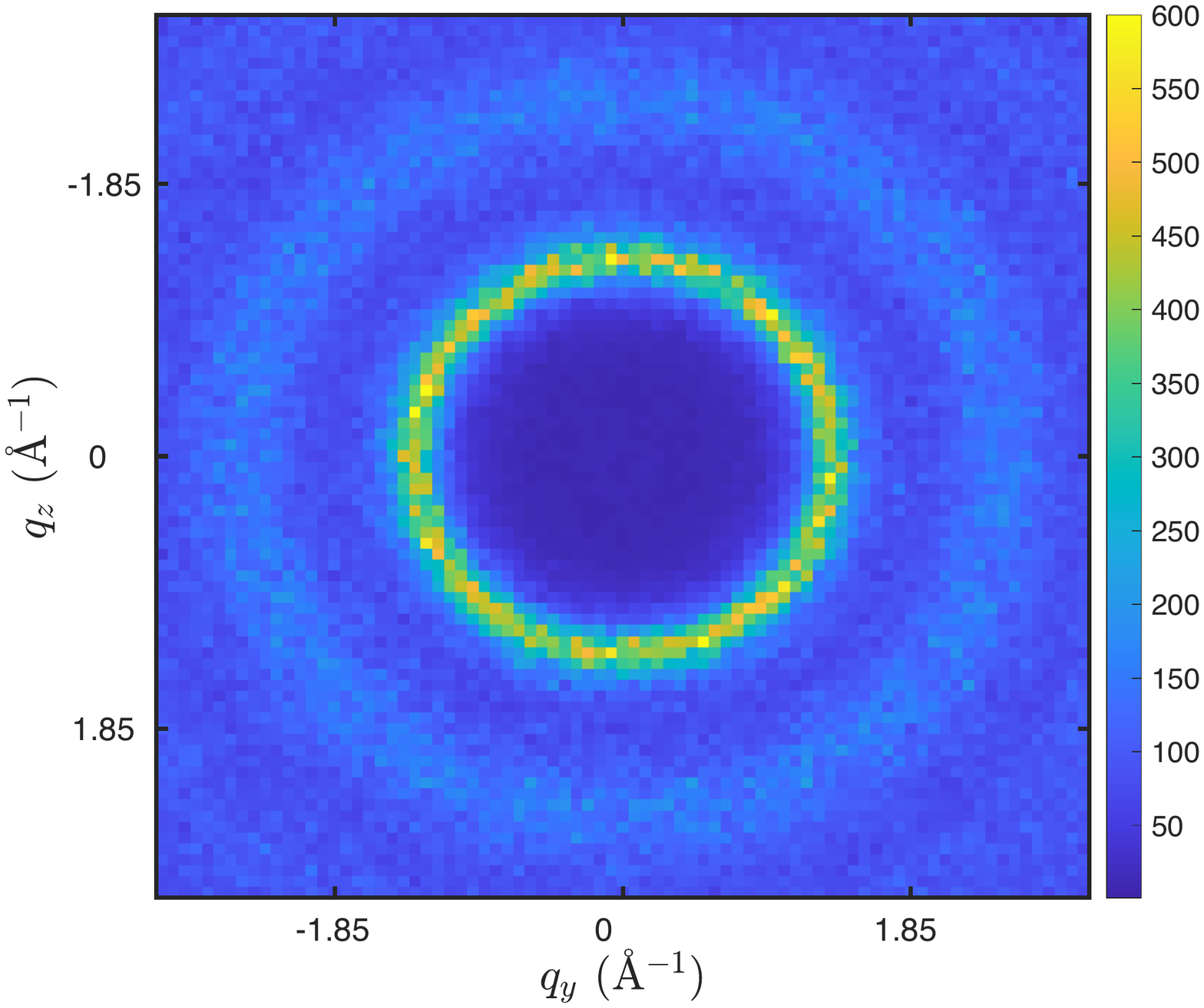}\label{fig:sqB}}
    \caption{The speckle patterns of scattered intensity $I(\bm{q})$ on a spherical slice in the $\bm{q}$-space (to keep $|\bm{k}_{f}| = |\bm{k}_{i}|$) for (a) the final configuration  and (b) time-averaged pattern over 40 snapshots with  $43.12$ ps separation between consecutive snapshots.
    \label{2dxray}}
\end{figure}

Fig.~\ref{2dxray} shows the computed scattered intensity on the detector grid from Algorithm 2. For the intensity computation for all atoms in the simulation box, the single snapshot shows the speckle pattern expected in an XPCS experiment (Fig.~\ref{fig:sqA}). The computation of the intensity averaged over $40$ frames (Fig.~\ref{fig:sqB}) starts to show the diffraction ring that we would expect to see at larger exposure times. In our analysis, we compute and store the intensity speckles from the instantaneous snapshots, similar to the data acquisition and storage procedure in XPCS experiments. 
The interval between consecutive snapshots is $43.12$ ps. 

\begin{table}[H]
    \centering
    \begin{tabular}{|l|c|}
    \hline
       \textbf{Method}  &  \shortstack{CPU Time (s) \\ per frame}   \\
       \hline
        Algorithm 1 (single point) & 0.0015 \\
        Algorithm 1 ($81\times81$ slice) & 4.1836 \\
        Algorithm 1 ($81\times 81\times 81 $ grid) & 416.8243  \\
        Algorithm 2 ($N_{\rm grid} = 400$)& 0.4527 \\
        \hline
    \end{tabular}
    \caption{CPU time taken by the parallelized implementation (over 32 CPUs) of the direct method (Algorithm 1) and FFT-based method (Algorithm 2), over 32 atomic configurations.}
    \label{tab:my_label}
\end{table}

The FFT-based method summarized by {\it Algorithm 2} has the advantage over the direct approach that the scattering intensity can be computed over the entire 3D grid all at once.
To demonstrate the computational efficiency of the FFT-based approach, we compare the CPU time taken by Algorithm 1 and Algorithm 2 in Table~\ref{tab:my_label}. In addition to the evaluation of X-ray intensity at a single $\bm{q}$ point, we also test the evaluation at a $81 \times 81$ 2D slice and a $81 \times 81 \times 81$ 3D grid in $\bm{q}$-space. 
The computation over the 2D slice is relevant because it contains the region of interest in the diffraction pattern, as shown in Fig.~\ref{fig:sqA}.
The computation over the 3D grid is also relevant because it can be used to quickly construct diffraction patterns corresponding to incident X-rays from all directions (e.g. to improve statistics of the computational predictions).
As mentioned Algorithm 1 computes the intensity over each scattering wavevector independently. Here the computation over the $\bm{q}$-space is parallelized on 32 CPUs. 
On the other hand, in Algorithm 2 a single CPU is used to compute the intensity over the entire 3D grid in $\bm{q}$-space for a given atomic configuration (frame).
Here the computation is parallelized on 32 CPUs over atomic configurations (snapshots).
Table~\ref{tab:my_label} shows that Algorithm 2 is much more efficient than Algorithm 1
when the scattering intensity on a large number of $\bm{q}$-points need to be evaluated for each frame.
Computing the statistical properties of the intensity speckles described in Section~\ref{subsec:stat_prop}, requires the computation of the intensity over a large number of frames ($5000$ frames in Sections~\ref{subsec:time_corr} and \ref{subsec:opt_con}) and over multiple grid points in the $\bm{q}$-space, making the FFT-based Algorithm 2 the preferred choice for computational modelling of XPCS experiments.

\subsection{Time correlation of XPCS signal} \label{subsec:time_corr}
\begin{figure}[H]
    \centering
    \subfigure[]{\includegraphics[width=0.48\textwidth]{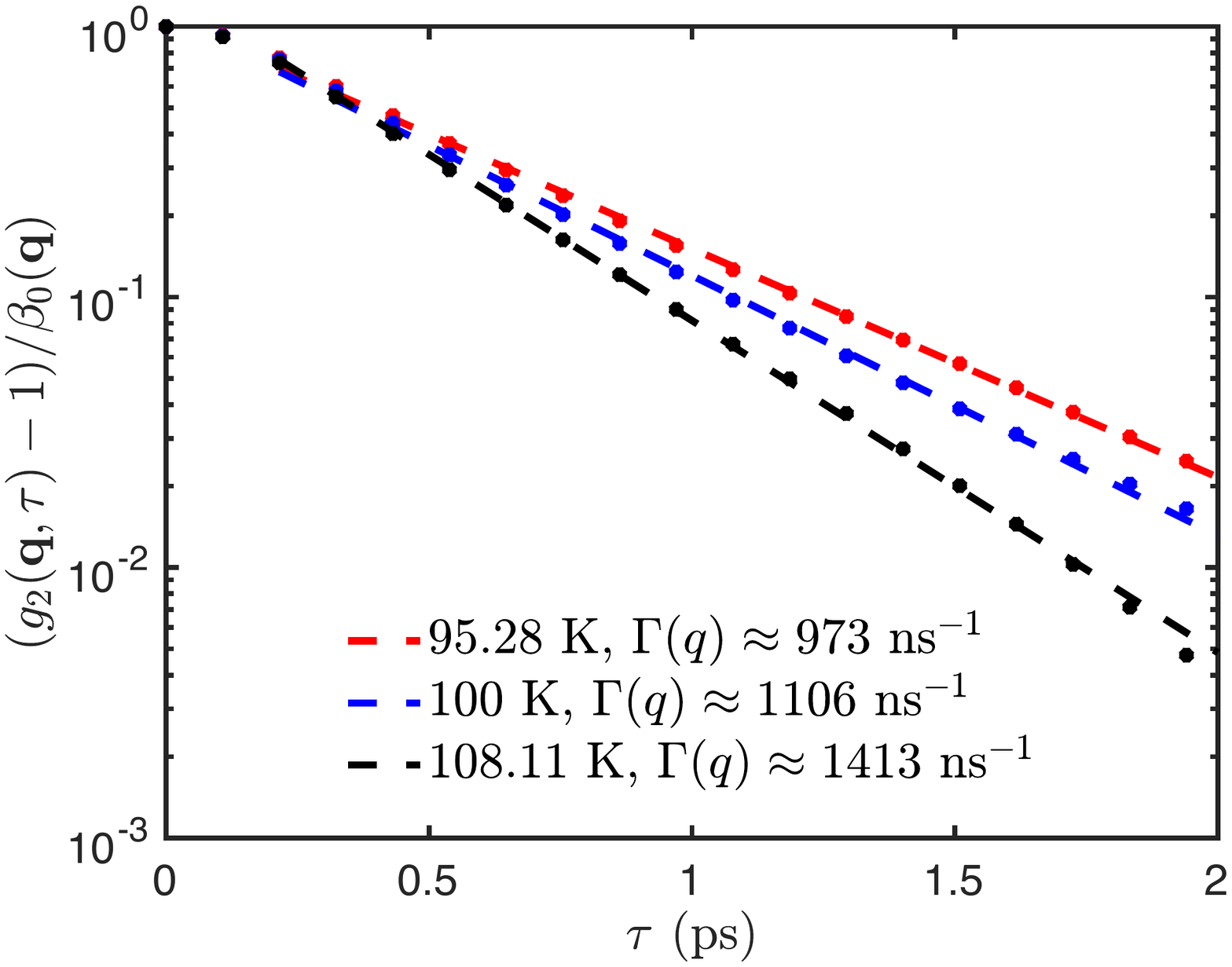}
    \label{fig:g2_all_45}}
   \subfigure[]{\includegraphics[width=0.48\textwidth]{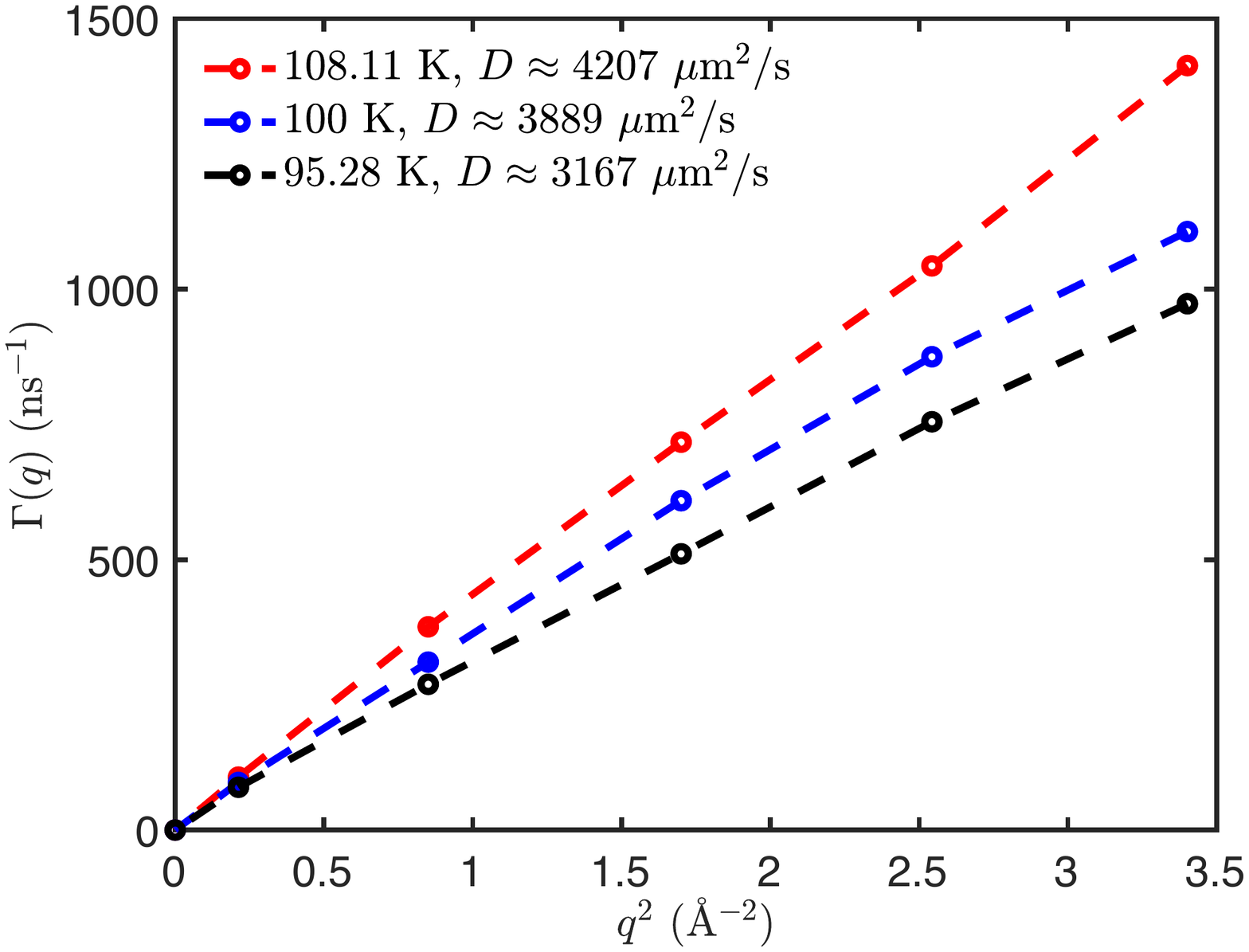}
    \label{fig:gamma_q2_45}}
    \caption{ (a) $(g_{2}(\bm{q},\tau)-1)/\beta_0(\bm{q})$ for $\lvert \bm{q} \rvert = 1.844 \pm 0.029$ \AA$^{-1}$ at $95.28$ K, $100$ K, and $108.11$ K. (b) $\Gamma(q)$ as a function of $q^{2}$ at $95.28$ K, $100$ K and $108.11$ K. 
    \label{fig:MSD_g_2_45}}
\end{figure}

In this section, we examine whether the time correlation of the computational XPCS speckles for liquid Ar decays exponentially, and whether the corresponding decay rates are related to the diffusion constant (see Section~\ref{subsec:stat_prop}).
\begin{figure}[H]
    \centering
    \subfigure[]{\includegraphics[width=0.48\textwidth]{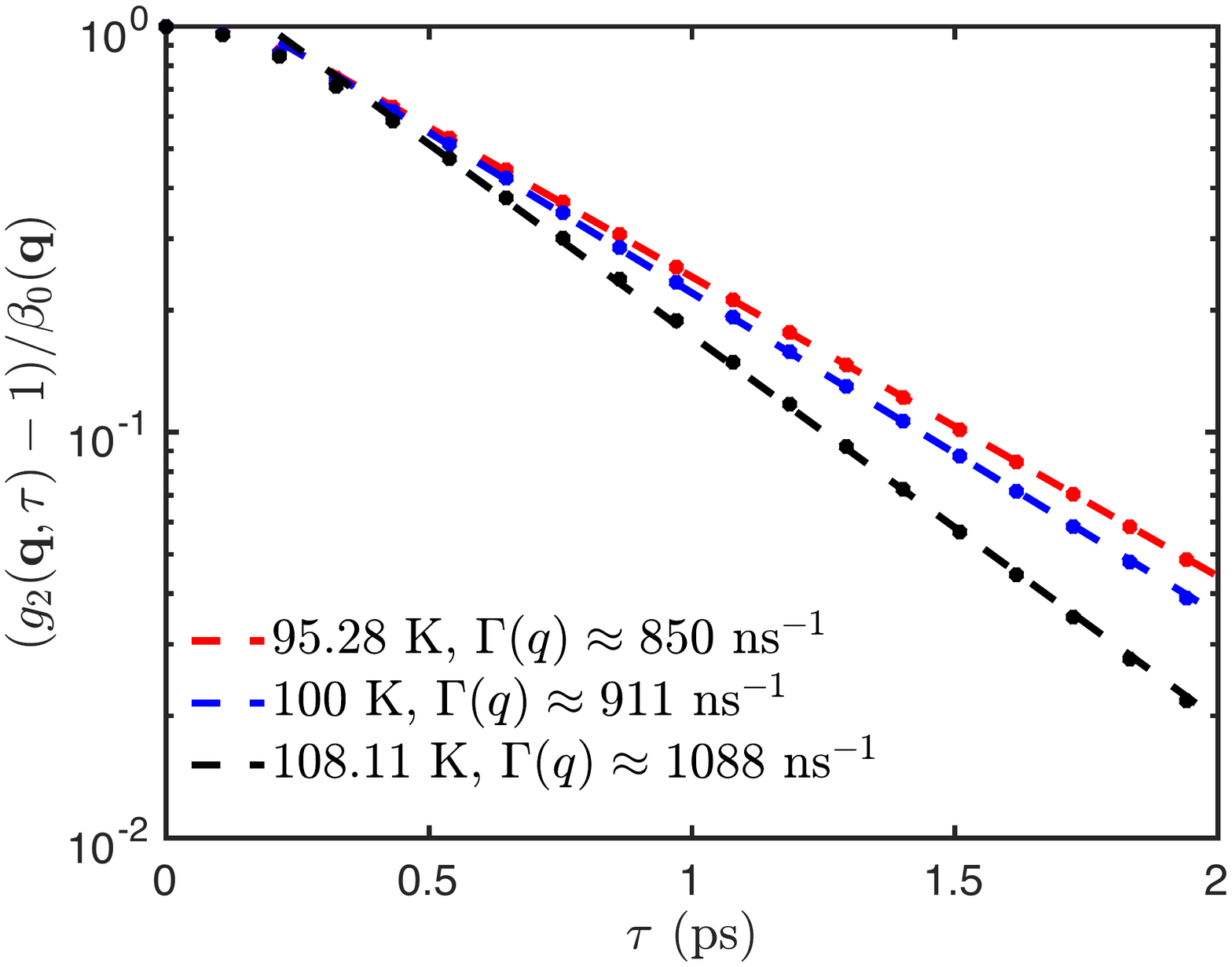}
    \label{fig:g2_all}}
   \subfigure[]{\includegraphics[width=0.48\textwidth]{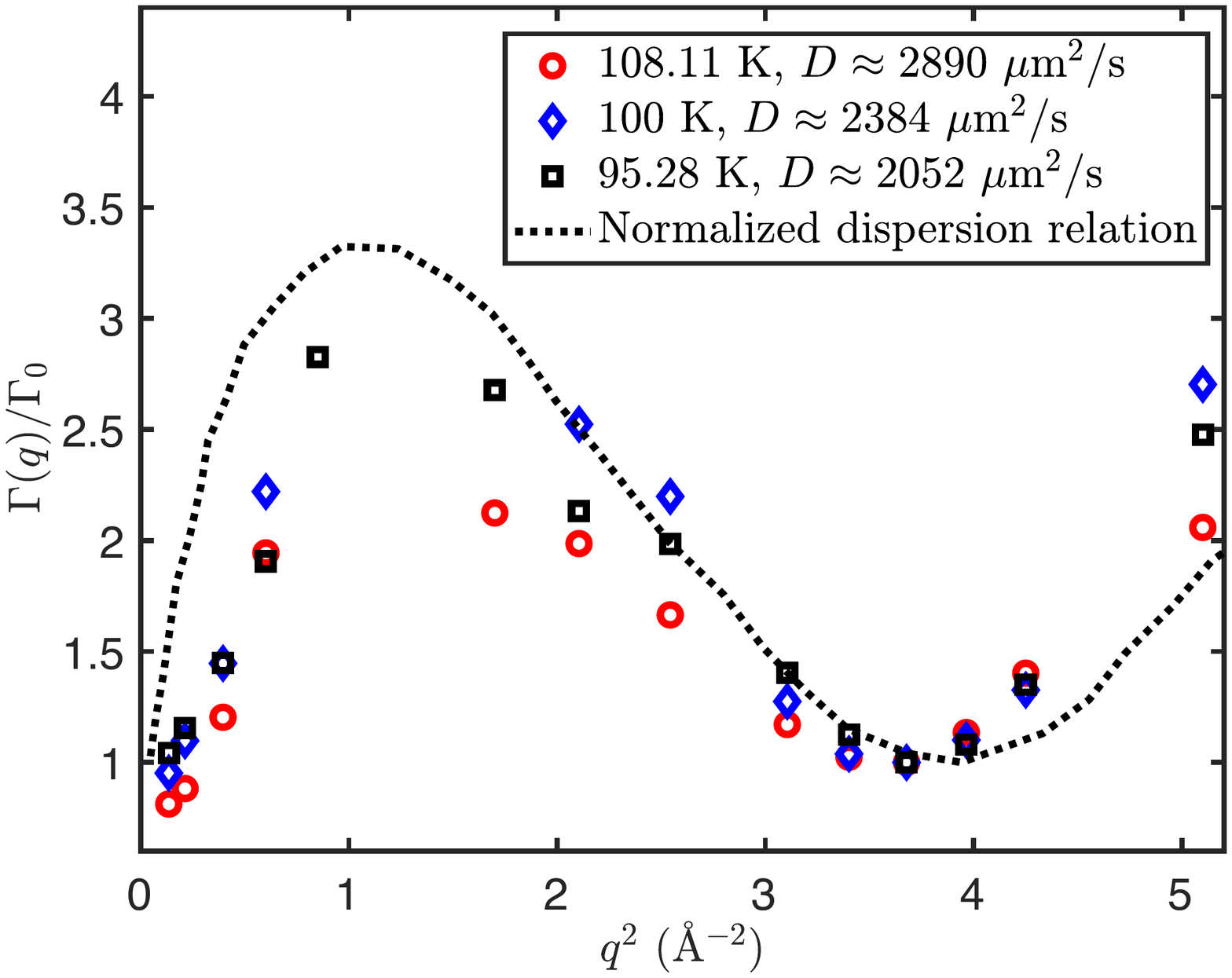}
    \label{fig:gamma_q2}}
    \caption{ (a) $(g_{2}(\bm{q},\tau)-1)/\beta_0(\bm{q})$ for $q = 1.844 \pm 0.029$ \AA$^{-1}$ at $95.28$ K, $100$ K, and $108.11$ K. (b) $\Gamma(q)\, /\Gamma_0$  as a function of $q^{2}$ at $95.28$ K, $100$ K and $108.11$ K, and the normalized analytical dispersion relation for liquid Ar~\cite{singwi1970collective}. $\Gamma_0$ is the value at the first minima which occurs at $q\approx 1.92$ \AA.
    \label{fig:MSD_g_2}}
\end{figure}
Fig.~\ref{fig:g2_all} shows that correlation function indeed decays exponentially with time for large $\tau$ for $q = 1.844 \pm 0.029$ \AA$^{-1}$, which corresponds to the wavevector magnitude around the first diffraction ring.
By fitting $(g_{2}(\bm{q},\tau)-1)/\beta_0(\bm{q})$ to $\exp[-\Gamma(\bm{q})\tau]$ in the large $\tau$ limit, we can extract the decay rates $\Gamma(\bm{q})$.
Fig.~\ref{fig:gamma_q2} plots the resulting $\Gamma(\bm{q})$ as a function of $q^2$; it can be seen that  $\Gamma(\bm{q})$  clearly deviates from the $Dq^2$ behavior expected from a purely diffusional process.
The behavior of $\Gamma(q)$, especially having a minimum around $q = 1.92$ \AA, is consistent with previous theoretical estimates and computational results on the intermediate scattering function~\cite{barker1972density,singwi1970collective}, where a non-linear dispersion relation is observed for liquid Ar at low temperatures. When we normalize the dispersion relation and $\Gamma(q)$ by their value at the local minimum around $q\approx 1.92$\,\AA, we see that data at all temperatures exhibit a qualitatively similar trend. 
It is interesting to note that, if we take the local minimum value of the $\Gamma(q)$ function, and divide by the corresponding $q^2$ value, we can obtain a rough estimate of the diffusivity $D$ that is within $40\%$ of the true diffusivity value computed from mean-square-displacements (MSD), as shown in Table~\ref{tab:MSD} (row for all atoms).
Therefore, our results suggest that through the $\Gamma(q)$ function, XPCS can provide valuable information on the intermediate scattering function, and even an order of magnitude estimate of the diffusivity in simple liquids.
\begin{figure}[H]
    \centering
    \subfigure[]{\includegraphics[width=0.44\textwidth]{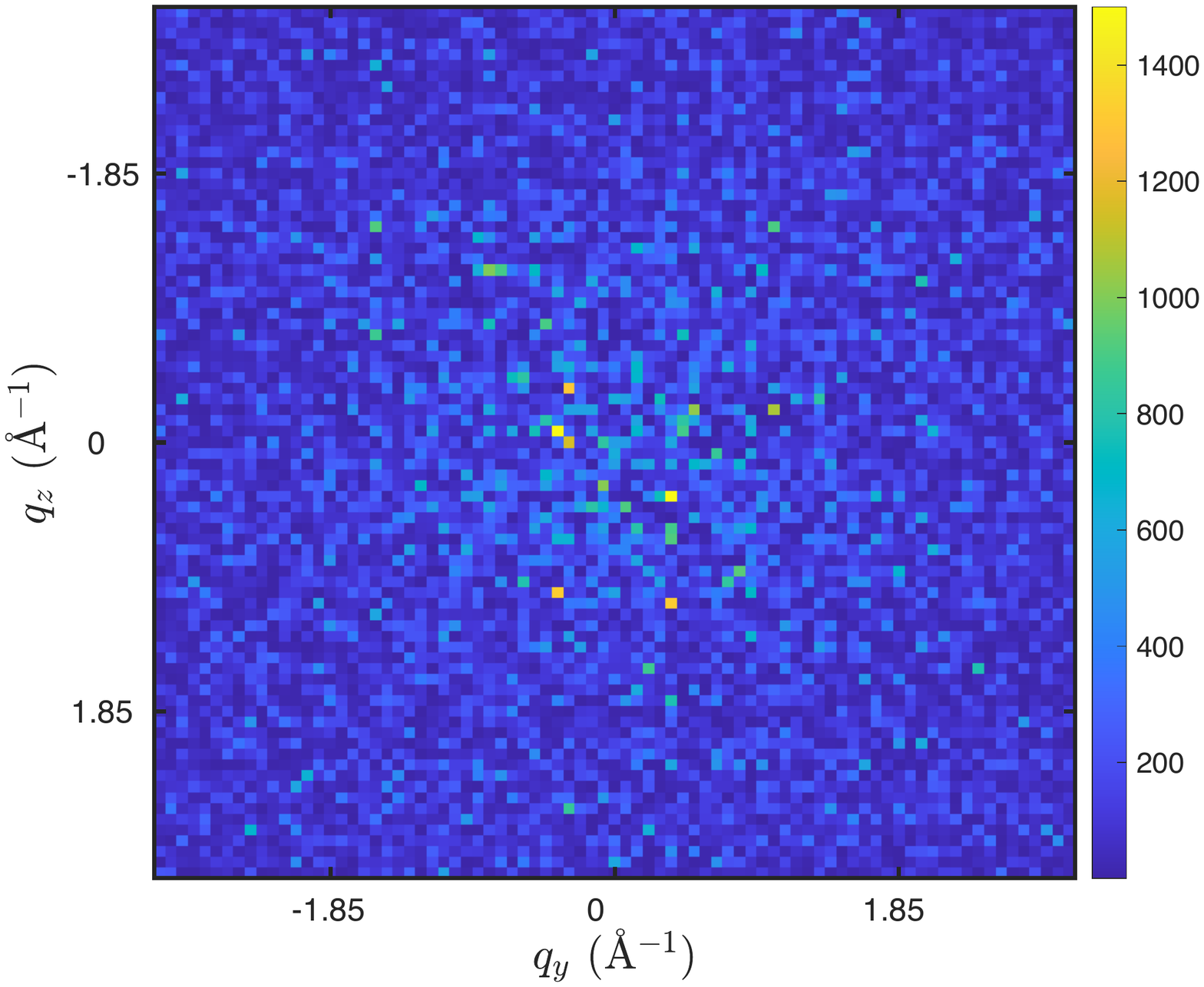}\label{fig:sqA_45}}
    \subfigure[]{\includegraphics[width=0.44\textwidth]{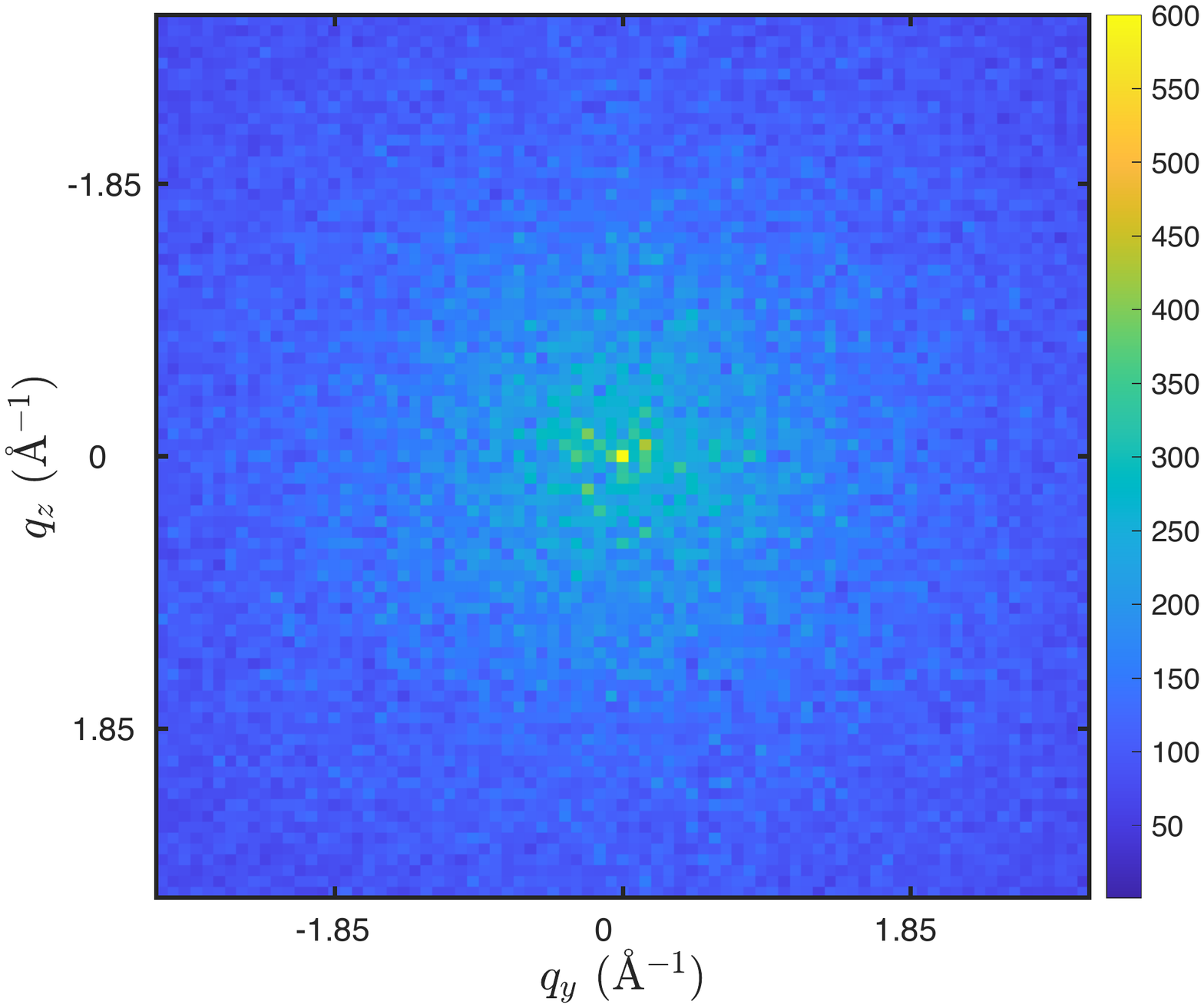}\label{fig:sqB_45}}
    \caption{The speckle patterns of scattered intensity $I(\bm{q})$ on a spherical slice in the $\bm{q}$-space (to keep $|\bm{k}_{f}| = |\bm{k}_{i}|$) from $45$ tracer atoms for (a) the final configuration  and (b) time-averaged pattern over 40 snapshots with  $43.12$ ps separation between consecutive snapshots.\label{2dxray_appendix}}
\end{figure}
To construct a benchmark case and a proof-of-concept in which $\Gamma(\bm{q})$ truly follows the $Dq^2$ trend in the low $q$ limit, we consider a thought experiment of tracer diffusion.
We randomly pick 45 atoms (as isotopes) out of the 4000 Ar atoms in the MD simulation cell.
For simplicity, we assume that these 45 atoms scatters X-ray much more strongly than the remaining atoms, whose scattering effect will thus be ignored.
The X-ray scattering speckle patterns due to these 45 atoms are shown in Fig.~\ref{2dxray_appendix}, which no longer shows the diffraction ring and looks very different from Fig.~\ref{2dxray}.
Fig.~\ref{fig:g2_all_45} shows that the correlation function for the resulting speckle patterns also decays exponentially with time.
Fig.~\ref{fig:gamma_q2_45} plots the corresponding decaying rate $\Gamma(\bm{q})$ as a function of $q^2$.
It can be seen that in this case $\Gamma(\bm{q})$ is indeed proportional to $q^2$ in the low $q$ limit.
Here we obtain estimates of the diffusion constant $D$ by fitting the $\Gamma(\bm{q})$ data over the entire $q$ range to $Dq^2 + D_2q^4$, in order to account for nonlinear effects at large $q$.
The results are shown in Table.~\ref{tab:MSD}, together with the $D$ values accurately computed from MD simulations based on the mean-square displacement (MSD).
\begin{table}[H]
    \centering
    \begin{tabular}{|l|c|c|c|}
    \hline
        \textbf{Method} & \shortstack{$D$ at $95.28$ K \\($\mu$m$^2$/s)} & \shortstack{$D$ at $100$ K \\ ($\mu$m$^2$/s) }& \shortstack{$D$ at $108.11$ K\\ ($\mu$m$^2$/s)} \\
        \hline
      MSD  & 2903& 3586& 4575\\
      XPCS (all atoms)  & 2052& 2384 & 2890\\
    XPCS (45 tracer atoms) & 3167& 3889& 4207\\
    \hline
    \end{tabular}
    \caption{The diffusivity, $D$, obtained from the MD trajectory using the mean-squared displacement (MSD) method and the computational XPCS method.}
    \label{tab:MSD}
\end{table}
Table.~\ref{tab:MSD} shows that in the hypothetical case of only a small fraction of tracer atoms dominate the X-ray scattering signal, the XPCS decay rates $\Gamma(q)$ can indeed provide accurate information of the diffusivity (within $8\%$).

\subsection{Optical contrast of XSVS signal} \label{subsec:opt_con}
\begin{figure}[H]
    \centering
    \subfigure[]{\includegraphics[width=0.38\textwidth]{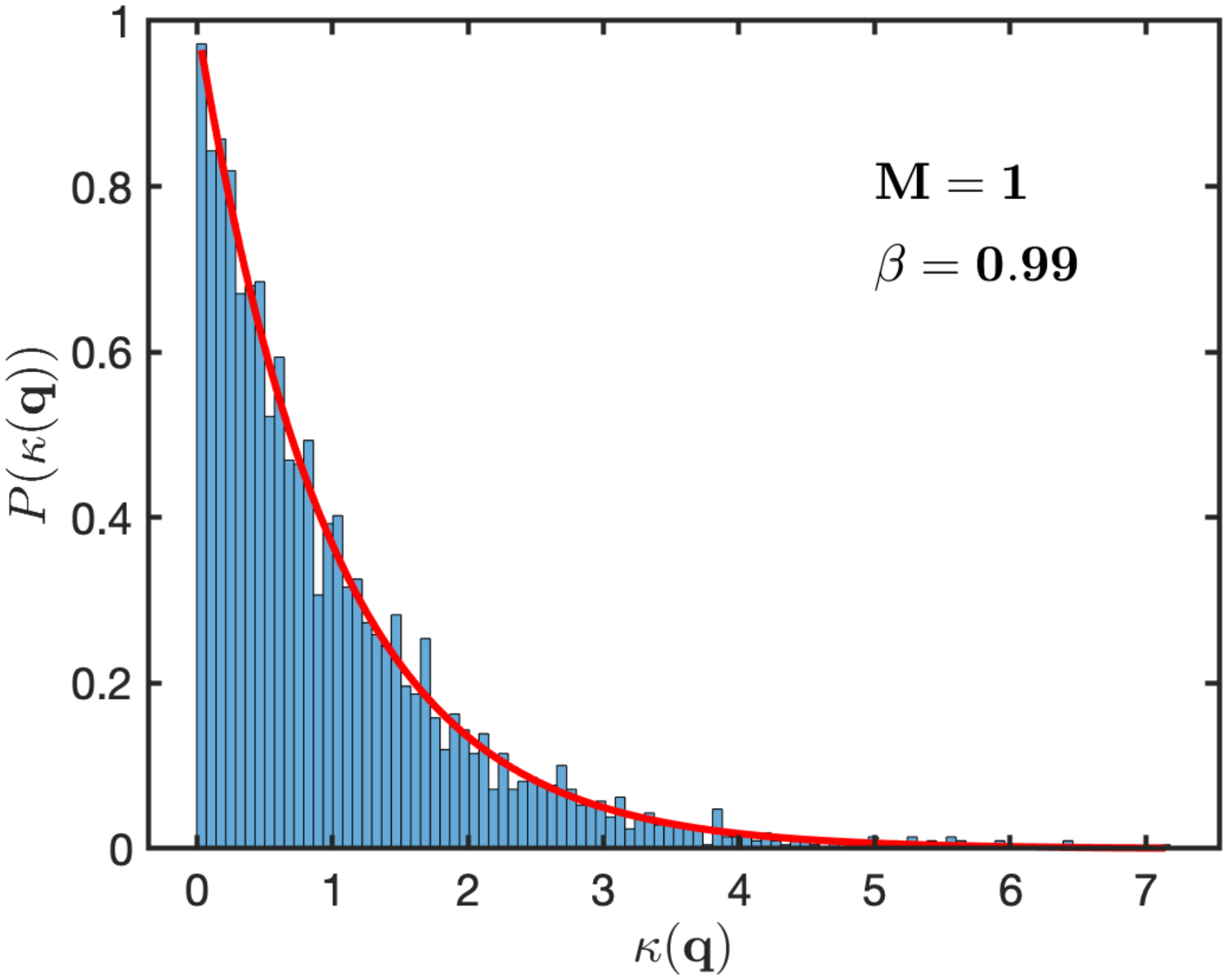}
    \label{fig:M_1}}
    \subfigure[]{\includegraphics[width=0.38\textwidth]{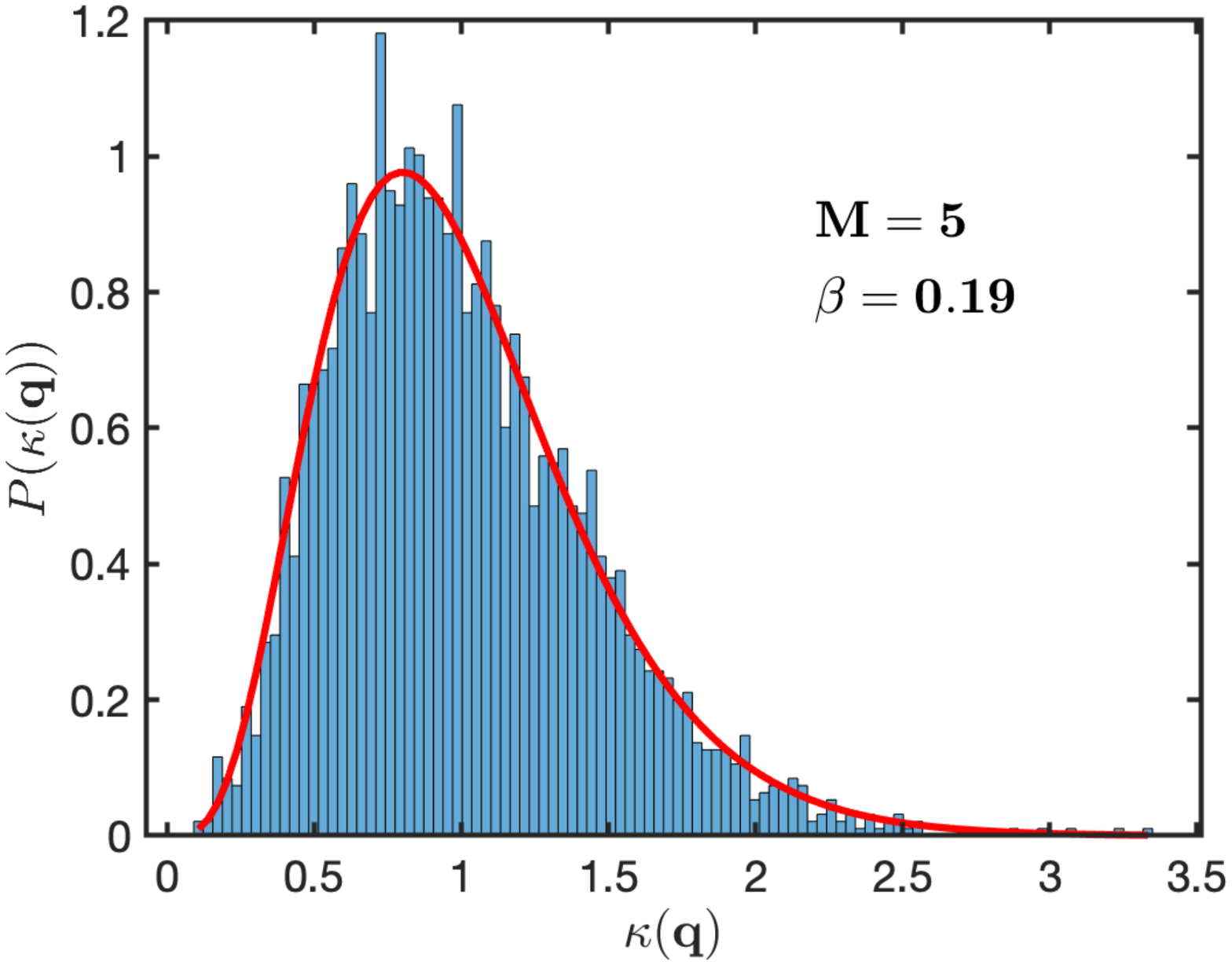}
    \label{fig:M_5}}
    
    \subfigure[]{\includegraphics[width=0.38\textwidth]{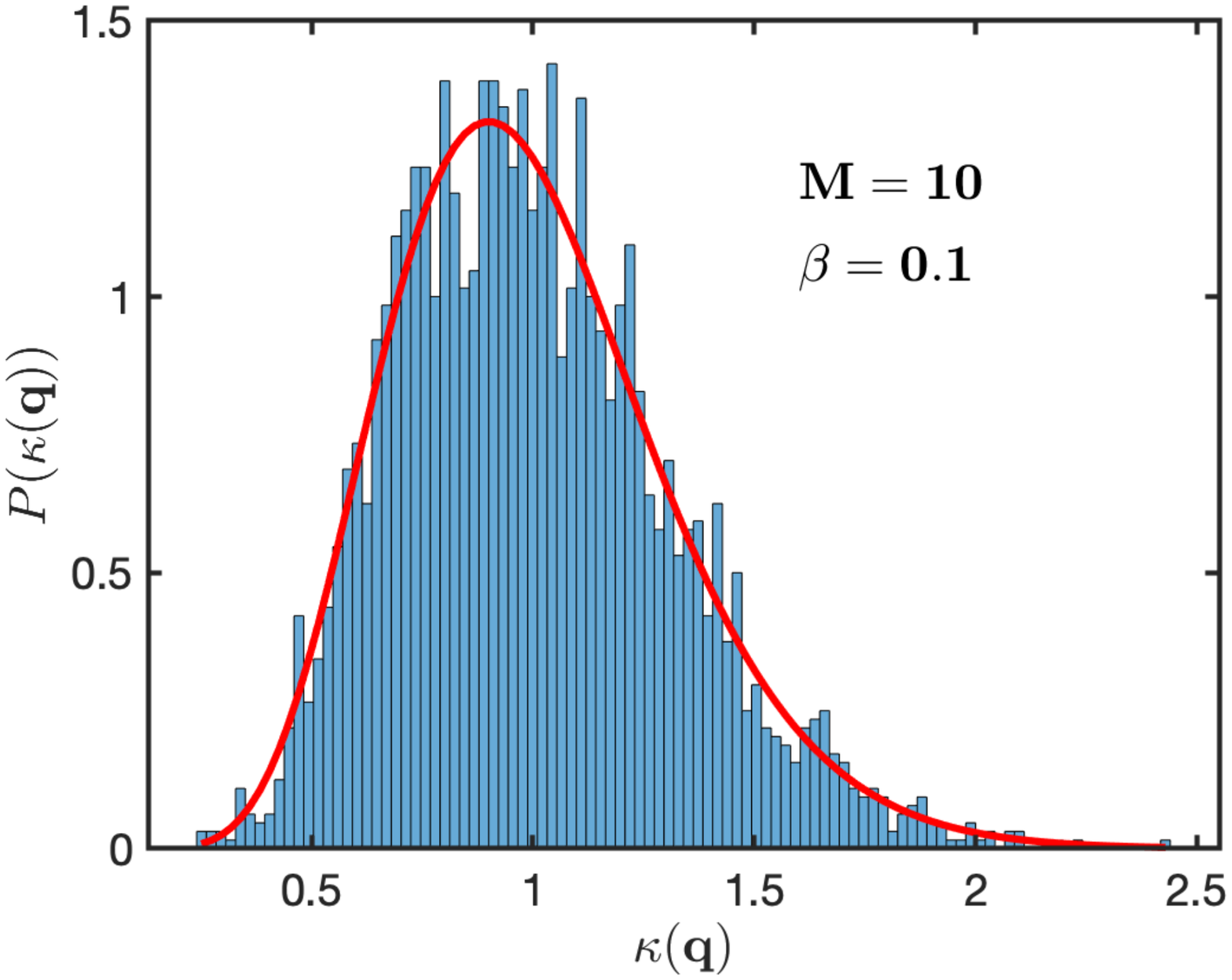}\label{fig:M_10}}
    \subfigure[]{\includegraphics[width=0.38\textwidth]{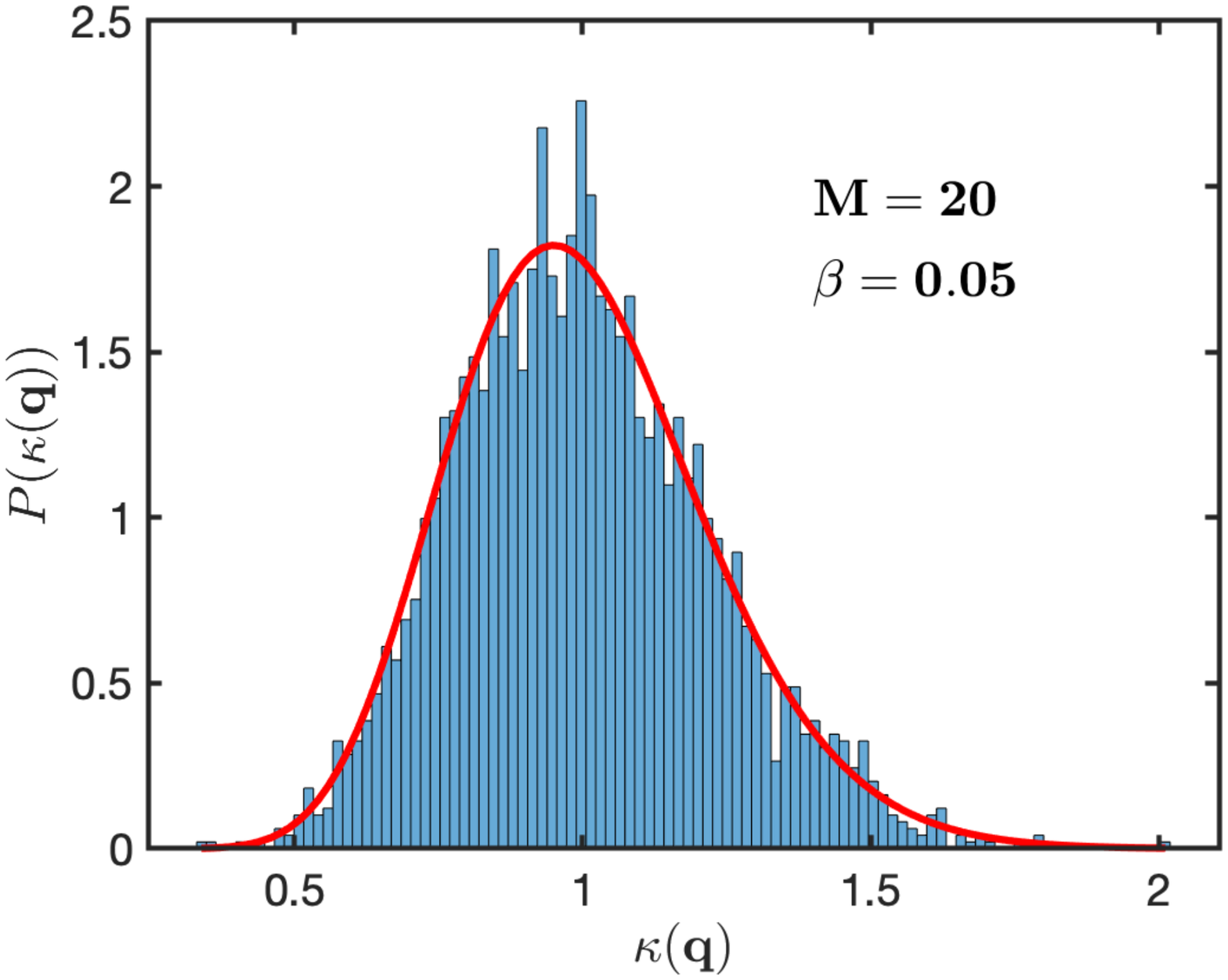}\label{fig:M_20}}
    \caption{The distribution of average normalized intensity for (a) $M=1$, (b) $M=5$, (c) $M=10$, and (d) $M=20$ independent configurations for $\lvert \bm{q} \rvert = 1.844 \pm 0.029$ \AA$^{-1}$. The red solid curve is the fit predicted by Eq.~(\ref{eq:Erlang}).  \label{fig:opt_cont}}
\end{figure}

In this section, we examine the statistical properties of the scattered intensities in a single computed speckle pattern, corresponding to a single incident X-ray pulse, as is done in XSVS experiments.
For an isotropic system, we expect the intensities at all $\bm{q}$-vectors of the same magnitude $q$ to have identical statistical distribution.
Fig.~\ref{fig:opt_cont}(a) shows the histogram of X-ray intensities for all $\bm{q}$-vectors whose magnitudes are in the range of $1.844 \pm 0.029$ \AA$^{-1}$.
It can be seen that the scattered intensities follow the exponential distribution,
\begin{equation}
    P\left(\kappa(\bm{q})\right) = {\exp }{\left[-\kappa(\bm{q}) \right]},
\end{equation}
where $\kappa(\bm{q}) = {I(\bm{q})}/{\langle I(\bm{q}) \rangle_{{q}}}$ and $P\left(\kappa(\bm{q})\right)$ is the probability distribution of $\kappa(\bm{q})$.
The exponential distribution is expected for perfectly coherent incident X-ray beam~\cite{dainty1975stellar} (as is assumed in the computational model).
The speckle contrast under a fully coherent beam is $\beta(q)=1$. This condition holds true as long as the sample volume is smaller than the coherence volume of the X-ray~\cite{zardecki1973volume} and the incident X-ray is perfectly monochromatic. 

\begin{figure}[H]
    \centering
        \subfigure[]{\includegraphics[width=0.48\textwidth]{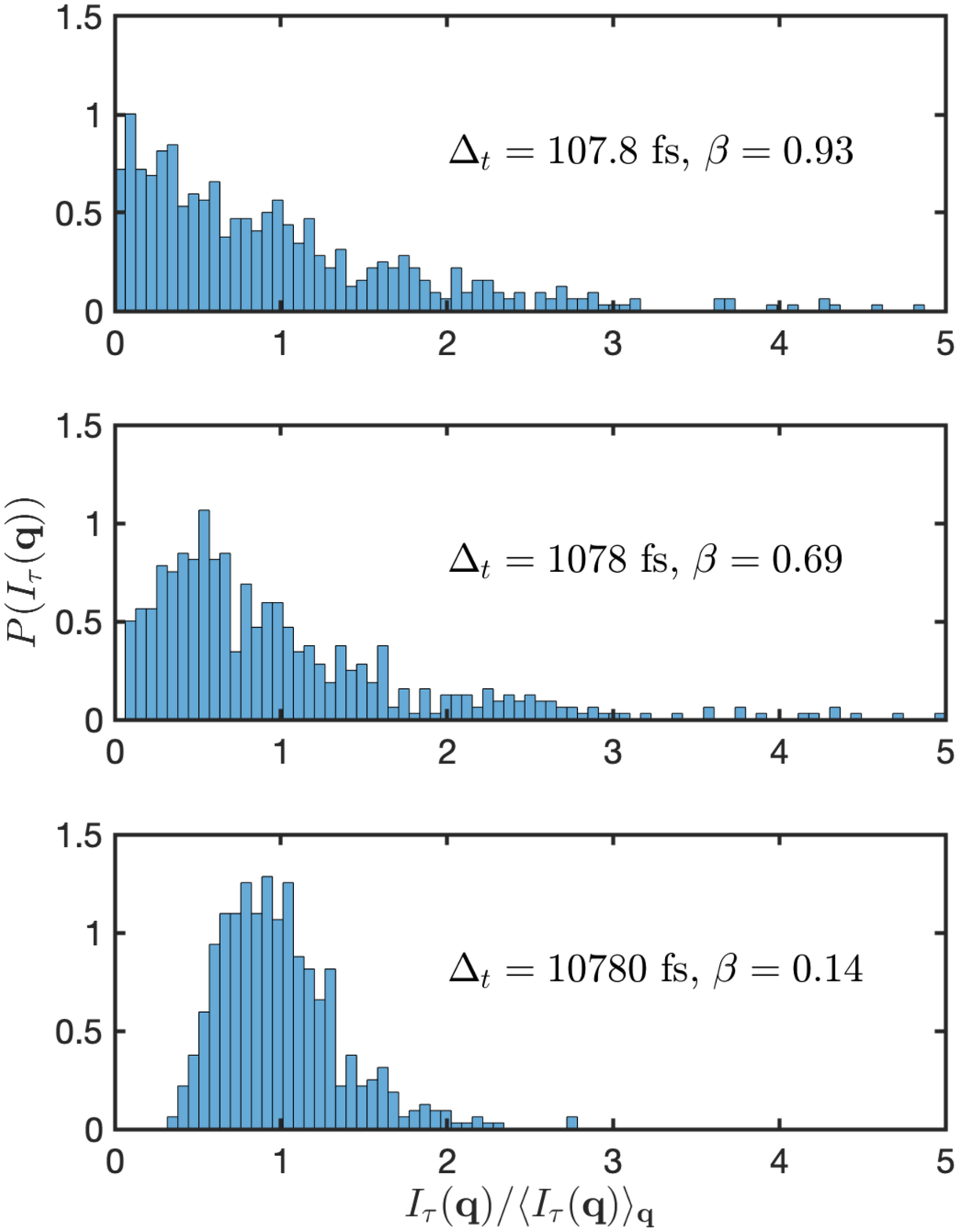}
    \label{fig:M_all2}}
    \subfigure[]{\includegraphics[width=0.48\textwidth]{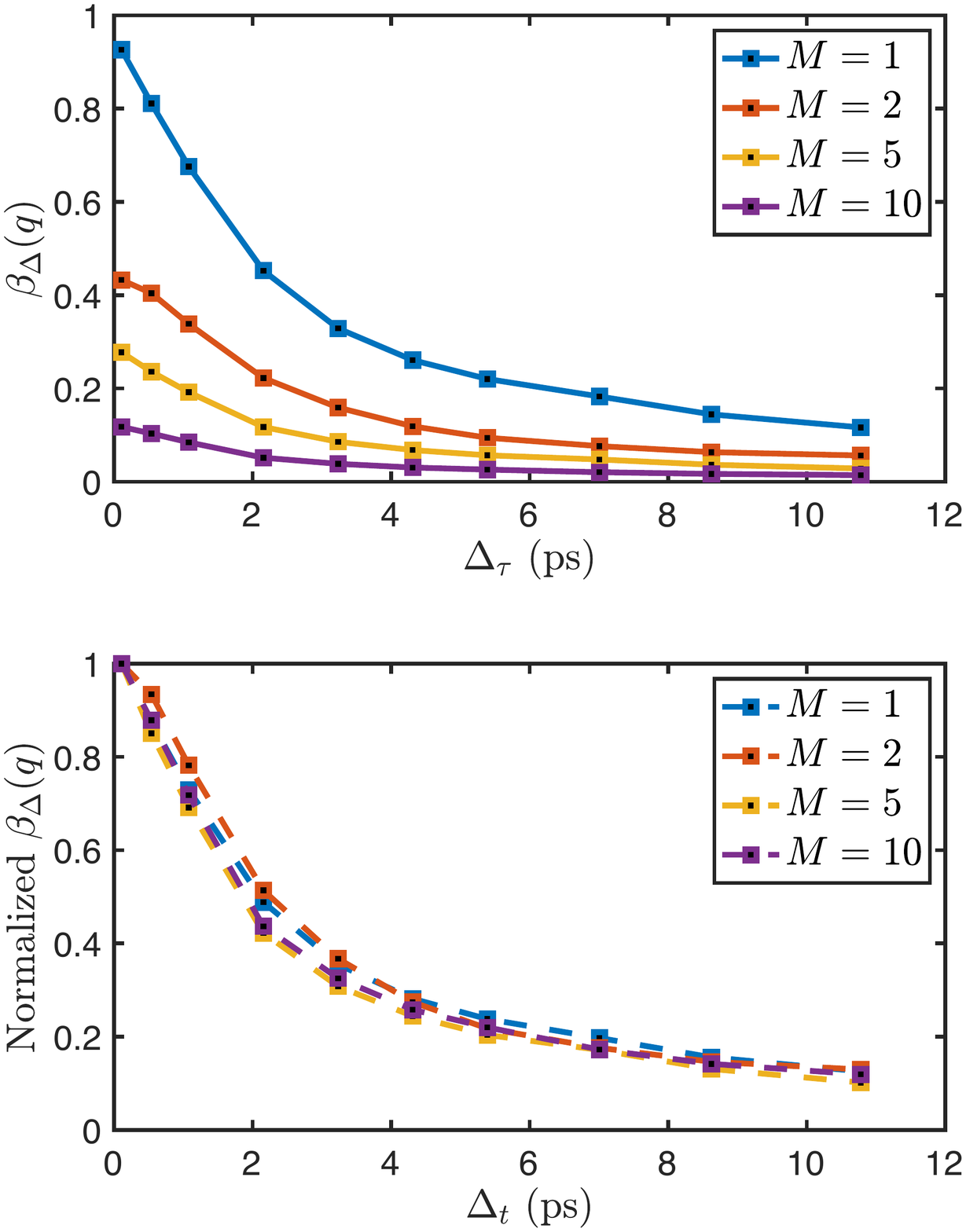}
    \label{fig:temp_all2}}
    \caption{(a) The distribution of the normalized integrated intensity for exposure times $\Delta_{t}=107.8,\,1078\, \text{ and}\,10780$ fs at $100$ K for one MD trajectory ($M=1$). (b) Optical contrast $\beta_{\Delta}(q)$ and ${\beta_{\Delta}(q)}/{\beta_{0}(q)}$ at different exposure times for a different number ($M$) of independent speckles at $100$ K. Both these calculations are carried out for $\lvert \bm{q} \rvert = 1.844 \pm 0.029$ \AA$^{-1}$.\label{fig:opt_cont_converge}}
\end{figure}

In practice, the X-ray beam used in the XPCS experiments may have a coherence volume smaller than the sample volume.
The observed speckle pattern would then be the incoherent superpositions (i.e.\,sum over intensity instead of complex amplitude) from several independent scattering volumes.
In computational XPCS, we can model this behavior by superimposing the scattering intensities from several independent atomistic configurations.
The resulting X-ray intensity is expected to satisfy the Erlang distribution (distribution of the sum of independent exponential random variables),
\begin{equation}
\label{eq:Erlang}
    P_{M}\left(\kappa(\bm{q})\right) =\frac{M^M\kappa(\bm{q})^{M-1}}{ \tilde{\Gamma}(M)} \,{\exp}{\left[-\kappa (\bm{q})M)\right]}, 
\end{equation}
where $\tilde{\Gamma}(M) = (M-1)!$ is the gamma function, and $P_{M}\left(\kappa(\bm{q})\right)$ is the probability distribution of the $\kappa(\bm{q})$ over $M$ independent speckle patterns.
The resulting optical contrast is expected to become $\beta(q) = 1/M$.
Fig.~\ref{fig:opt_cont}(c)-(d) show the histograms of scattered intensity computed by superimpositing $M = 5, 10, 20$ independent atomistic configurations.
These atomistic configurations are extracted from the same MD simulation but are separated by multiples of time duration $\tau_{\rm sample}=54$ ps, which is much greater than the correlation time $\tau_{\rm cor}$ of the Ar liquid.
It can be seen that the distributions of the scattered X-ray intensity agree well with the expected Erlang distribution; the corresponding optical contrast $\beta(q)$ also agrees well with the expected value of $1/M$.
Next, we examine the effect of finite duration $\Delta_t$ of the incident X-ray pulse on the optical contrast.
Fig.~\ref{fig:M_all2} show the histogram of the scattered intensity in the $q$-range $\lvert \bm{q} \rvert = 1.844 \pm 0.029$ \AA$^{-1}$, for X-ray pulse duration of $\Delta_t = 107.8$, 1078, 10780\,fs (assuming perfectly coherent X-ray beam, i.e. $M = 1$).
Fig.~\ref{fig:temp_all2} shows the optical contrast $\beta(q)$ as a function of pulse duration $\Delta_t$.
As expected, $\beta(q)$ decreases with increasing $\Delta_t$; the time at which the optical contrast decreases to half of its maximum value gives an estimate of the correlation time, which is about $2$~ps here.
Fig.~\ref{fig:M_all2} also shows the combined effect of limited coherence volume (i.e. $M$ > 1) and finite pulse duration.
It can be seen that for all $M$ values, the optical contrast $\beta(q)$ decreases with pulse duration $\Delta_t$.
When normalized by the optical contrast at the $\Delta_t \to 0$ limit, the decay of optical contrast with $\Delta_t$ for all cases collapse onto a single master curve, as is commonly assumed when analyzing XSVS data.
\subsection{Computational XSVS of water}
In this section, we examine the efficacy of the FFT-based model in studying other liquid systems by verifying the decay of the optical contrast of water with increasing pulse duration~\cite{perakis2018coherent}. 
For the simulation of water, we use the TIP4P/2005  and follow the simulation procedure described by \citet{perakis2018coherent}. The simulations are run at seven different temperatures ($250$ K, $290$ K, $296$ K, $328$ K and $330$ K) to compare with previous work~\cite{perakis2018coherent}.
\begin{figure}[H]
    \centering
    \subfigure[]{\includegraphics[width=0.48\textwidth]{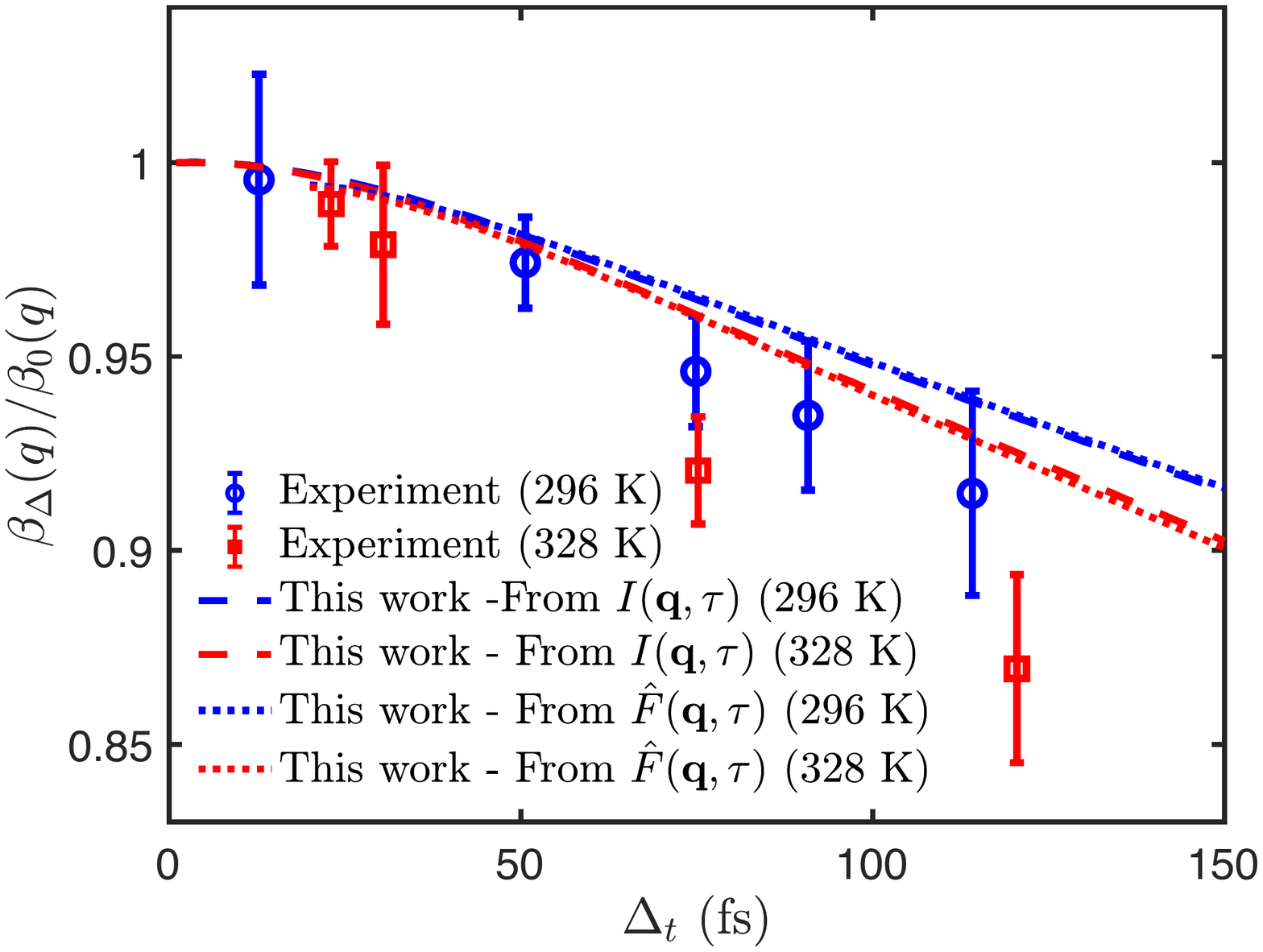}\label{fig:beta_water}} 
    \subfigure[]{\includegraphics[width=0.48\textwidth]{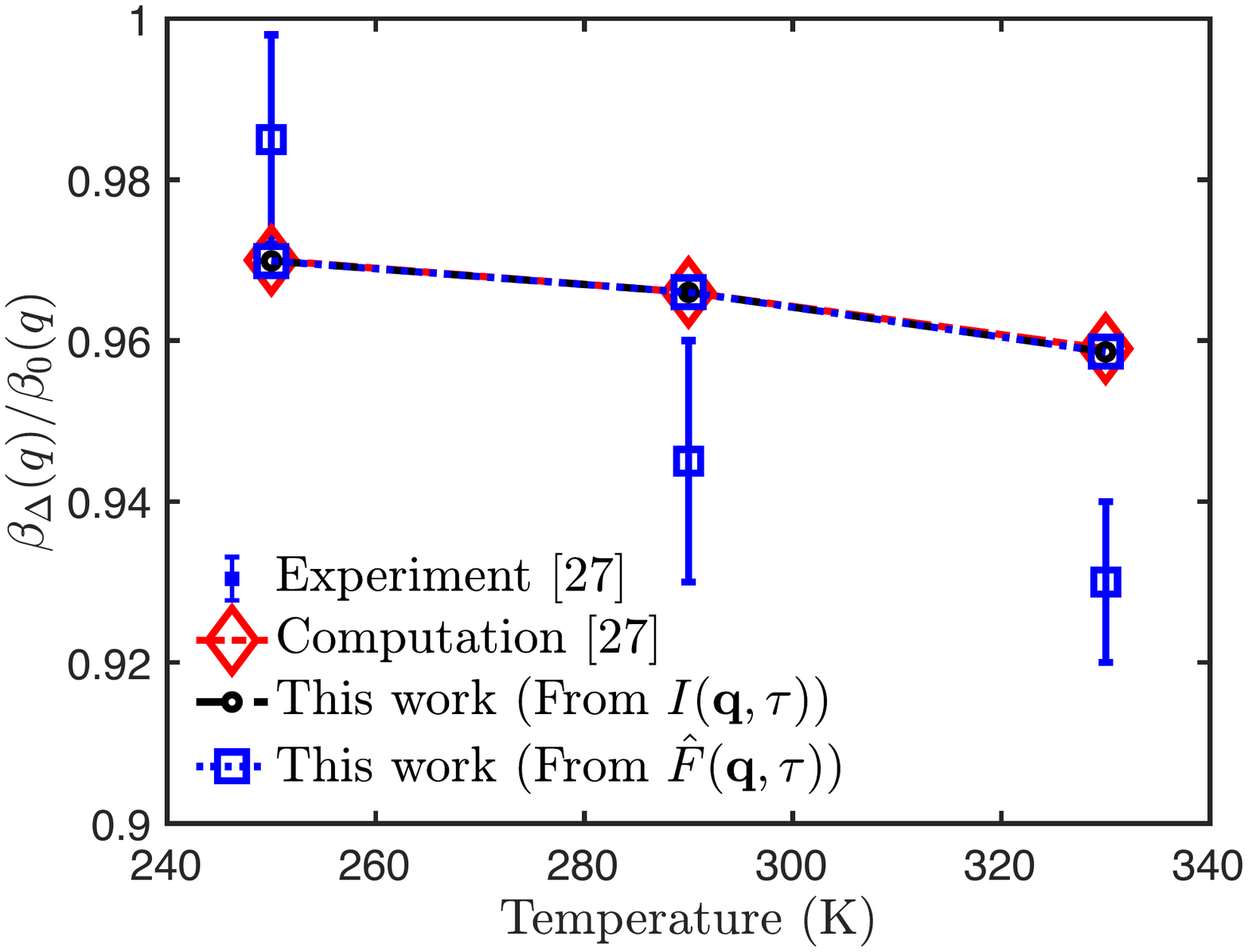}\label{fig:beta_temp}} 
    \caption{(a) The relative optical contrast, $\beta_{\Delta}(q)/\beta_{0}(q)$ as a function of the exposure time, $\Delta_{t}$ observed experimentally~\cite{perakis2018coherent} (red squares and blue circles) and from our computational model (dashed lines) [computed using Eqs.~(\ref{eq:bq_Iq}) and (\ref{eq:beta_alt_int})]. These calculations are carried out for $\lvert \bm{q} \rvert = 1.95 \pm 0.1$ \AA$^{-1}$. (b) Comparison of the computationally obtained $\beta_{\Delta}(q)/\beta_{0}(q)$ at $\Delta_{t}=75$ fs against the experimental results and computational results reported by ~\citet{perakis2018coherent}.}
    \label{fig:beta_water_temp}
\end{figure}
We compute the $\beta_{\Delta}(q)/\beta_{0}(q)$ for simulations at $296$ K and $328$ K at $\lvert \bm{q} \rvert = 1.95 \pm 0.1$ \AA$^{-1}$ by using Eqs.~(\ref{eq:bq_Iq}) and (\ref{eq:beta_alt_int}), where Eq.~(\ref{eq:bq_Iq}) computes the $\beta_{\Delta}(q)$ from the $I(\bm{q},t)$ and Eq.~(\ref{eq:beta_alt_int}) computes the $\beta_{\Delta}(q)$ from the $\hat{F}(\bm{q},t)$. 
Our results are in close agreement with the computational predictions based on the intermediate scattering function $\hat{F}(\bm{q},t)$~\cite{perakis2018coherent} and the Siegert relation, as shown in Fig.~\ref{fig:beta_water} (the experimental data from \cite{perakis2018coherent} are shown by markers).
On the other hand, the FFT-based method does not assume or make use of the Siegert relation. 
Instead, we compute the optical contrast $\beta_{\Delta}(q)$ directly from $I(\bm{q},t)$, in the same way as in the experiments.

At higher temperatures, we expect the correlation time to decrease due to the increase in thermal fluctuations. The decrease in correlation time results in a quicker decay in the optical contrast $\beta_{\Delta}(q)$, as shown in Fig.~\ref{fig:beta_water}. 
The difference between the experimental data and computational predictions is likely due to the inaccuracy of the interatomic potential for water, as noted earlier~\cite{perakis2018coherent}.
Our results demonstrate that the optical contrast can be computed directly from the scattering intensity $I(\bm{q})$, like in the XSVS experiments, instead of through the intermediate scattering function $F(\bm{q}, t)$, under the assumption of the validity of the Siegert relation for the system under analysis.
Additionally, the computation using Eq.~(\ref{eq:bq_Iq}) can be generalized to other pulse shapes of the incident X-ray (not restricted to a square wave).

\section{Conclusions} 
\label{sec:conclusion}
This paper compares two methods for computing XPCS and XSVS signals from molecular dynamics (MD) configurations, and shows that the computational signals satisfy the expected relations and properties of XPCS and XSVS experiments. We demonstrated the equivalence of the FFT-based method with the direct method. The FFT-based method has higher efficiency due to the simultaneous computation of the intensity signal over the entire grid in $\bm{q}$-space. 
We provide numerical evidence for the Siegert relation for the XPCS signals computed from MD simulations of liquid Ar. We also show the equivalence between the optical contrast $\beta(q)$ computed from the $\bm{q}$-space and from the time domain. 

Through the computational results, we show that the time correlation of XPCS speckle intensities can provide information on dynamical properties such as dispersion relation and diffusivity.
When applied to the XSVS experiments, the computed speckle pattern exhibits optical contrast of unity in the limit of infinitesimally short X-ray pulses and a decay of optical contrast with increasing exposure time. 
We show that using the FFT-based method one can obtain the same numerical result on the XSVS measurement for water without invoking the Siegert relation as previous reports using a different method that requires the Siegert relation to hold.

Through these numerical examples, we confirm the equivalence and validity of the direct method and the FFT-based method for computational XPCS and XSVS. 
The model systems under study are liquid Ar and water, whose dynamics are probed by molecular dynamics (MD) simulations at the femtosecond to nanosecond timescales. 
The approach presented here can be extended to coarse-grained molecular dynamics (CGMD) simulations, which can help the intepretation of XPCS signals of polymeric networks at the microseconds to milliseconds timescales in terms of the bond and segmental dynamics of polymer chains.
%
%
 
\section*{Conflict of Interest}
The authors declare no competing interests.
\section*{Data Access}
All data is available in the main text and the Supplementary appendices. Further information about the computation can be obtained on request to the corresponding author.
\section*{Acknowledgements}
Use of the Linac Coherent Light Source (LCLS), SLAC National Accelerator Laboratory, is supported by the U.S. Department of Energy, Office of Science, Office of Basic Energy Sciences under Contract No. DE-AC02-76SF00515. C.B.C. acknowledges support from the Department of Defense (DoD) through the National Defense Science and Engineering Graduate (NDSEG) Fellowship Program.
S.M. and W.C. acknowledge support from the Precourt Pioneering Project of Stanford University.

\end{document}